\documentclass[12pt]{amsart}
\usepackage{amsmath,amsthm}
\input amssym.def
\input amssym

\usepackage{hyperref}
\hypersetup{
    pdftitle={Memory Estimation of Inverse Operators},    
    pdfauthor={Anatoly Baskakov}{Ilya Krishtal},     
    pdfsubject={Memory Estimation of Inverse Operators},   
    pdfcreator={Ilya Krishtal},   
    pdfproducer={Ilya Krishtal}, 
    pdfkeywords={Banach modules}{Beurling spectrum}{Wiener's lemma}, 
     colorlinks,%
    citecolor=blue,%
    filecolor=black,%
    linkcolor=magenta,%
    urlcolor=black
}

\chardef\bslash=`\\ 

\makeatletter
\def\verbatim{\interlinepenalty\@M \@verbatim
   \leftskip\@totalleftmargin\advance\leftskip2pc
   \frenchspacing\@vobeyspaces \@xverbatim}
\makeatother
\newtheorem{thm}{Theorem}[section]
\newtheorem{cor}[thm]{Corollary}
\newtheorem{lem}[thm]{Lemma}
\newtheorem{prop}[thm]{Proposition}

\newtheorem{ass}[thm]{Assumption}

\theoremstyle{definition}
\newtheorem{defn}{Definition}[section]

\theoremstyle{remark}
\newtheorem{rem}{Remark}[section]
\newtheorem{exmp}{Example}[section]

\numberwithin{equation}{section}

\newcommand{\begeq}{\begin {equation}}
\newcommand{\eq}{\end{equation}}
\newcommand{\bs}{\begin {split}}
\newcommand{\es}{\end{split}}
\newcommand{\bp}{\begin {prop}}
\newcommand{\ep}{\end {prop}}
\newcommand{\bt}{\begin {thm}}
\newcommand{\et}{\end {thm}}
\newcommand{\bc}{\begin {cor}}
\newcommand{\ec}{\end {cor}}
\newcommand{\bl}{\begin {lem}}
\newcommand{\el}{\end {lem}}
\newcommand{\bpf}{\begin {proof}}
\newcommand{\epf}{\end {proof}}
\newcommand{\bi}{\begin {itemize}}
\newcommand{\ei}{\end {itemize}}
\newcommand{\ben}{\begin {enumerate}}
\newcommand{\een}{\end {enumerate}}
\newcommand{\brem}{\begin {rem}}
\newcommand{\erem}{\end {rem}}

\newcommand{\bd}{\begin {defn}}
\newcommand{\ed}{\end {defn}}

\newcommand{\bex}{\begin {exmp}}
\newcommand{\eex}{\end {exmp}}
\newcommand{\bai}{b.a.i.~}
\newcommand{\norm}[1]{\left\|{#1} \right\|}

\newcommand{\fhat}{\hat{f}}
\newcommand{\eps}{\epsilon}

\newcommand{\cB}{\mathfrak B}

\newcommand{\A}{\mathcal A}

\newcommand{\C}{\mathcal C}

\newcommand{\B}{{\mathcal B}}
\newcommand{\F}{\mathcal{F}}

\newcommand{\cP}{{\mathcal P}}
\newcommand{\cQ}{{\mathcal Q}}

\newcommand{\Y}{\mathcal{Y}}
\newcommand{\TTT}{{\T\kern-.44em \T}}
\newcommand{\tTTT}{\widetilde{\T\kern-.44em \T}}
\newcommand{\chT}{\check\T}

\newcommand{\ZZ}{{\mathbb Z}}
\newcommand{\TT}{{\mathbb T}}
\newcommand{\RR}{{\mathbb R}}
\newcommand{\CC}{{\mathbb C}}
\newcommand{\NN}{{\mathbb N}}

\newcommand{\X}{\mathcal{X}}
\newcommand{\Z}{\mathcal{Z}}

\newcommand{\GG}{\mathbb{G}}

\newcommand{\G}{\mathbb G }

\newcommand{\HH}{{\mathcal H}}

\renewcommand{\d}{\delta}
\newcommand{\s}{\sigma}
\renewcommand{\l}{\lambda}
\renewcommand{\L}{\Lambda}
\renewcommand{\a}{\alpha}

\newcommand{\lnhat}{\widehat{L_{\nu}}}

\newcommand{\lnloc}{\widehat{L_{\nu}}^{loc}}

\newcommand{\g}{\gamma}
\newcommand{\Gg}{\widehat{\G}}
\newcommand{\T}{\mathcal{T}}
\newcommand{\W}{\mathcal{W}}

\newcommand{\ee}{{\textbf{h}}_{\alpha}}
\newcommand{\eet}{{\textbf{e}}_\alpha }
\newcommand{\bh}{{\textbf{h}}}

\DeclareMathOperator{\supp}{supp}

\newcommand{\aee}{\mbox{\rm \ae}}

\begin{document}

\title {Memory estimation of inverse operators}

\author{Anatoly G. Baskakov and Ilya A. Krishtal}
\address{Department of Applied Mathematics and Mechanics, Voronezh State University, Russia 394693 \\
email: mmio@amm.vsu.ru}
\address{Department of Mathematical Sciences, Northern Illinois University, DeKalb, IL 60115 \\
email: krishtal@math.niu.edu}

\thanks{ The first author is supported in part by RFBR grant 10-01-00276.
The second author is supported in part by NSF grant DMS-0908239.}


\date{\today }




\begin{abstract}
We use methods of harmonic analysis and group representation theory to estimate memory decay of the inverse operators in Banach spaces. The memory of the operators is defined using the notion of the Beurling spectrum. We obtain a general continuous non-commutative version of the celebrated Wiener's Tauberian lemma with estimates of the ``Fourier coefficients'' of inverse operators. In particular, we generalize various estimates of the elements of the inverse matrices. The results are illustrated with a variety of examples including integral and integro-differential operators.
\end{abstract}
\maketitle

\section{Introduction}
Wiener's Tauberian Lemma  \cite{W32} is a classical result in harmonic analysis which states that if a periodic function $f$ has an absolutely convergent Fourier series and never vanishes
then the function $1/f$ also has an absolutely convergent Fourier series. This result has many extensions (see \cite{ABK08, B08, BK10, B90, B97Sib, B97Izv, GKW89, GK10, GL06, J90, K11, K90, L53, S95, S78, S05, S07} and references therein), 
which have been used, for example, in such diverse areas as differential equations \cite{K99}, pseudo-differential operators \cite{GR08, GS07},
 frame theory \cite{ABK08, BCHL06I, F09, G04, KO08, S07, S08ACM}, time-frequency analysis \cite{GL04, GL06, KO08}, sampling theory \cite{AAK09, ABK08, SS09, S10}, finite-section method \cite{GRS10, RRS98}, etc.  

A standard reformulation of Wiener's lemma states that if an invertible operator is defined by a (bi-infinite) Laurent matrix with
summable diagonals then the inverse has the same property. In \cite{B90, GKW89, K90} it was shown (independently) that the matrix
does not have to be Laurent. Hence, if one interprets a matrix entry $a_{ij}$ as a ``memory cell'', that is information on what impact an event at ``time'' $j$ has on the state of the system at ``time'' $i$, the Wiener property for an operator roughly means that such an operator does not spread the information too far or has a localized memory. In \cite{K11} one of us showed 
that Wiener's lemma 
is really a statement about the preservation of memory localization by inverse operators, and the way the memory is defined is
irrelevant to a large extent. In this paper we continue and expand this line of research paying particular attention to the case of ``continuous
memory'' and obtaining specific memory estimates for the inverse operators.

The following Wiener's lemma extension developed in \cite{B97Izv} serves as a starting point and an example for the research in this paper.

Let $\X, \Y$ be infinite dimensional complex Banach spaces and $\cP=(\cP_k)$ and $\cQ =(\cQ_k)$, $k\in S\subseteq\ZZ$, be two sequences of (continuous) idempotents from the Banach algebras $B(\X)$ and $B(\Y)$ of bounded linear operators on $\X$ and $\Y$, respectively. We assume that each of the sequences is disjunctive, i.e. $\cP_k\cP_j = \delta_{jk}\cP_k$, $k,j\in S$, and similarly for the sequence $\cQ_k$, $k\in S$. As usual, by $\delta_{jk}$ we mean the Kronecker delta. We also assume that for every $x\in \X$ and $y\in\Y$ we have
\[x = \sum_{k\in S} \cP_kx\quad\mbox{ and }\quad y = \sum_{k\in S} \cP_ky,\]
where the series converge unconditionally, and that 
\begeq\label{crp}
C(\cP) = \sup_{\a_k\in\TT}\norm{\sum_{k\in S} \a_k\cP_k}<\infty\mbox{ and } 
C(\cQ) = \sup_{\a_k\in\TT}\norm{\sum_{k\in S} \a_k\cQ_k}<\infty,\eq
where $\TT = \{\g\in\CC:\, |\g|=1\}$ is the unit circle.

A sequence of idempotents with the above properties will be called a \emph{disjunctive resolution of the identity}. Given such $\cP$ and $\cQ$ one can define a matrix $\A = (A_{mn})$ for each operator $A$ in the vector space $L(\X,\Y)$ of bounded linear operators from $\X$ to $\Y$. The matrix is a mapping
$\A: S\times S\to L(\X,\Y)$ given by
\[\A(m,n) =A_{mn} =  \cQ_mA\cP_n,\ m,n\in S,\]
and the operators $A_{mn}$, $m,n\in S$, are called \emph{operator blocks}. If $A$ is continuously invertible, then the operator blocks of the operator $B = A^{-1}\in L(\Y,\X)$ are defined by
$B_{mn} = \cP_mB\cQ_n$, $m,n\in S$.

In \cite{B97Izv} the class of operators with summable diagonals, which in this paper we will denote by $\W$, was defined using the sequence
\[d_A(k) = \sup_{m-n=k}\norm{A_{mn}},\ k\in S-S,\, \mbox{ and } d_A(k)=0,\mbox{ if } k\in\ZZ\backslash (S-S).\]
More precisely, we call
\[\W = \W(\X,\Y) = \{A\in L(\X,\Y):\, \sum_{k\in\ZZ}d_A(k)< \infty\}\]
the \emph{Wiener class} of operators (with respect to the resolutions of the identity $\cP$ and $\cQ$).

Wiener's lemma extension proved in \cite{B97Izv} states that if $A\in \W(\X,\Y)$ is 
invertible, that is
$B = A^{-1}\in L(\Y,\X)$, then $B\in\W(\Y,\X)$.
\brem
It was shown in \cite{K11} that if $\X=\Y=\HH$ is a Hilbert space, then the result remains valid if the resolutions of the identity $\cP$ and $\cQ$ are not assumed to be disjunctive.
\erem

Although we have avoided mentioning it until now, the key technique for estimating the norms of the operator blocks of the inverse operator in \cite{B97Izv} is based on the spectral properties of the representations of locally compact Abelian (LCA-) groups. 
Indeed, the resolutions of the identity $\cP$ and $\cQ$ give rise to two $2\pi$-periodic representations
\[U_\X:\RR\to B(\X) \mbox{ and } U_\Y:\RR\to B(\Y)\]
defined by
\begeq\label{exres}
U_\X(t) x = \sum_{k\in S} e^{ikt}P_kx,\  U_\Y(t) y = \sum_{k\in S} e^{ikt}Q_ky, \ x\in\X, y\in\Y, t\in\RR.
\eq
The above representations, in turn, produce another $2\pi$-periodic representation
$U_{\X\Y}: \RR\to B(L(\X,\Y))$ via
\begeq\label{typop1}
U_{\X\Y} (t)A= U_\Y(t)AU_\X(-t).
\eq The \emph{Fourier series} of the operator $A$ is then given by
\[U_{\X\Y}(t) A \simeq \sum_{k\in\ZZ} e^{ikt}A_k,\]
where the \emph{Fourier coefficients} satisfy \cite{B97Izv}
\[A_k = \frac1{2\pi}\int_0^{2\pi} e^{-ikt}U_{\X\Y}(t) A \,dt = \sum_{m-n = k} \cQ_mA\cP_{n},\ k\in\ZZ.\]
The above integral converges in the strong operator topology and the series converges strongly and unconditionally.
If $C(\cP) = C(\cQ) = 1$, which is always true for some equivalent renormalization of $\X$ and $\Y$, we have $d_k(A) = \|A_k\|$ and the connection between the Wiener's class $\W$ and summability of Fourier coefficients becomes apparent. 

Often it is natural to refer to the Fourier coefficients of an operator as its \emph{diagonals} and the set of indices of non-zero diagonals as \emph{memory}. Obviously both notions are dependent upon the choice of the representation.
Although the paper \cite{B97Izv} deals only with representations \eqref{exres}, it is easily seen that the technique works for any periodic representation (or, equivalently, a representation of the LCA-group $\TT$). 
In \cite{BK10} it is shown how to prove analogous results in the almost periodic case (for representations of the Bohr compactification $\RR_c$ of $\RR$). In this paper we are primarily interested in theory based on representations of $\RR^d$, $d\in\NN$. In this case, memory is no longer a discrete set, the representations in \eqref{exres} are usually not well defined, and we have to use far more sophisticated methods of the spectral theory
of Banach modules (representations of LCA-groups). Fortunately, most of the required methods
have already been developed in \cite{B04, BK05}.

The remainder of the paper is organized as follows. In section \ref{nota1} we summarize the
notions and results of the spectral theory of representations of LCA-groups. 
In Section \ref{Cla} we introduce different classes of vectors with spectral decay such as Wiener Class, Beurling Class and others. We believe that some of the classes (one-sided exponential spectral decay, Sobolev-type classes) have not been previously considered in conjunction with Wiener's lemma type results. We also don't believe that any of these classes have been introduced at this level of generality.
Section \ref{Clop} introduces standard \cite{B04, BK05, BR75} Banach module structures on the space of operators which allow us to represent memory decay of  linear operators as the spectral decay of vectors in Banach modules. In Sections \ref{Wcomp}--\ref{sobol} we prove results about the memory decay of inverses to operators in different classes. The main Theorems are \ref{1sidexpwin}, \ref{invexp}, \ref{maint}, \ref{absB}, and \ref{sobolt}. It is worth noting here that some of these results are quantitative. While we use a combination of traditional methods of proof
(holomorphic extensions, Neumann series ``boot-strap'', Brandenburg trick, etc. \cite{G10}), 
we also employ a closed operator functional calculus that has not been used before. More importantly,
the framework of Banach modules (group representations) makes our results extremely general as illustrated by  
a variety of examples in Section \ref{examp}. We show that many of the previously known results of this kind can be obtained using our technique. On the other hand, an absolute majority of these examples cannot be obtained at this level of generality in other frameworks known to us. In particular, the block-matrix technique of \cite{B97Izv} fails because
in general the representation in \eqref{exres} is not well defined (the series may not converge unconditionally). Other matrix-based techniques \cite{S05, S07, S10ca} do not apply essentially for the same reason. Known results for integral operators \cite{K99, K01, S08} apply to narrower classes of operators. Finally,
 the abstract technique based on the Bochner-Phillips theorem \cite{B97Sib, BP42} has either not been developed for this level of generality or
cannot be used directly, because we consider operators in $L(\X,\Y)$ as well as Banach algebras with a group of automorphisms.

\section{Banach modules over weighted group algebras and their spectral properties}\label{nota1}

In this section we introduce the notation and develop the necessary tools from representation theory. The presentation is largely
analogous to that of \cite{BK05}, where isometric representations are used. Our interest in differential operators, however, forces
the use of unbounded representations, which requires an extra effort. Although, most of the results
in Sections \ref{Wcomp}--\ref{sobol} are formulated for bounded representations we prefer to
present the prerequisite theory more generally in order to use it in the sequel to this paper.

As above, $\X$ will denote a complex Banach space and $B(\X)$ will be the Banach algebra of
bounded linear operators on $\X$. By $\G$ we will denote an LCA-group
and by $\Gg$ -- its Pontryagin dual, the group of continuous unitary characters of $\G$. 
We write the operation additively on all LCA-groups, except when the group $\TT = \{\theta\in\CC:\ |\theta| = 1\}$ is used.
Next, we introduce the following definition of weights.

\begin{defn} \label{wegt}
A \emph{weight} is an even function $\nu: \G\to [1,\infty)$ such that
\[ \nu(g_1+g_2)\leq \nu(g_1)\nu(g_2),\ \mbox{ for all }\ g_1,g_2\in\G. \]
A weight is \emph{admissible} if it satisfies the \emph{GRS-condition}
\[\lim\limits_{n\to\infty}n^{-1}\ln\nu(ng)=0,\quad \mbox{ for all }\quad g\in\G,\quad ng = \underbrace{g+g+\ldots+ g}_{n\, {\rm times}}. \]
A weight is \emph{exponential} if $\GG$ is a subgroup of $\RR^d$ (possibly, with discrete topology) and
$\nu(g) = e^{a|g|}$ for some $a > 0$. 
For 
$g = (g_1,g_2,\ldots,g_d)\in\RR^d$ we will use the notation $|g|_p = \left(\sum\limits_{k=1}^d |g_k|^p\right)^{\frac1p}$, $1\le p< \infty$,  $|g| = |g|_1$, and $|g|_\infty = \max\limits_{1\le k\le d} |g_k|$. 
\end{defn}
 
We denote by $L_\nu(\G)$ the Banach algebra of (equivalence classes of) complex functions
on $\G$ integrable with weight $\nu$ with respect to the Haar measure on $\G$. If $\nu \equiv 1$, we obtain the standard space $L^1(\G) \equiv L_1(\G)$. We shall use both notations interchangeably, hinting at the possibility of using more general weights when the subscript is used.  The role of
multiplication is played by the convolution of functions and the norm is given by
\[\norm{f}_\nu = \int_\G |f(g)|\nu(g)dg .\] 
We refer to \cite{K09, L85} for the various properties of the algebra $L_\nu(\G)$ used below.

In this paper we assume that the weight $\nu$ is non-quazi-analytic \cite{D56, LMF73, L85}, that is 
\[\sum_{n=0}^\infty \frac{\ln \nu(ng)}{1+n^2} < \infty \quad \mbox{for all } g\in\G.\]
This condition ensures \cite{D56, K09} that the group algebra $L_\nu(\G)$ is regular and its spectrum 
 is isomorphic to the dual group $\Gg$.

We denote by $\fhat : \Gg \to \CC$ the Fourier transform of $f \in L_\nu(\G)$. The inverse Fourier
transform of a function  $h: \Gg \to \CC$ is denoted by $\check h$ or $h^\vee$. If $\GG = \RR^d$, we use the Fourier transform in the form
\[\hat f(\l) = \int_{\RR^d} f(t)e^{-it\cdot\l} dt,\ \l\in\RR^d.\]

The memory of the operators will eventually be defined using the properties of group representations
$\T: \G \to B(\X)$. Together with the representation in \eqref{exres}, the following representations are used most often.

If $\X$ is an appropriate space of functions (or distributions) on $\G$, the \emph{translation} representation $\T = T: \G\to B(\X)$ is given by
\begeq\label{trans1}
T(t)x(s) = x(s+t),\ x\in\X, \ s,t\in\G,
\eq 
and the \emph{modulation} representation $\T = M: \Gg \to B(\X)$ is defined by
\begeq\label{mod1}
M(\xi)x(s) = s(\xi)x(s),\ x\in\X, \ \xi\in\Gg, \ s\in \G.
\eq 
Other important examples of representations used in the paper are introduced below.

We assume that $\X= (\X,\T)$ is a non-degenerate Banach $L_\nu(\G)$-module \cite{B04, BK05},
with the module structure associated with a representation $\T: \G \to B(\X)$.
In particular, in this paper we consider only modules for which the following assumption is satisfied.

\begin{ass}\label{assu}
The following three conditions hold for the Banach $L_\nu(\G)$-module $\X$:
\begin{enumerate}
	\item the module $\X$ is non-degenerate, that is, if $fx = 0$ for all $f \in L_\nu(\G)$ then $x = 0$;
  \item  the module structure on $\X$ is associated with a representation $\T: \G \to B(\X)$, that is,
  for all $f \in L_\nu(\G)$, $x \in \X$, and $g\in\G$,
\begeq\label{redef}
\T(g)(fx) = \left(T(g)f\right)x = f\left(\T(g)x\right),
\eq  
where  $T$ is the translation representation (\ref{trans1});
  \item 
  $\norm{fx} \leq \|f\|_\nu\|x\|$, for all $f \in L_\nu(\G)$ and $x \in \X$.
\end{enumerate}
\end{ass}

\brem\label{cnepr}
We observe that  to ensure that $x=0$ it is enough to show that $fx = 0$
for every $f\in L_\nu(\GG)$ with $\supp \hat f$ compact. Indeed, since every $h\in L_\nu(\GG)$ 
is a limit of a net $(f_\a)$ of functions with $\supp\hat f_\a$ compact \cite{K09, L85}, we would get
\[hx = \lim_{\a} (h*f_\a)x = 0\]
and, hence, $x = 0$ since $\X$ is non-degenerate. 
\erem

If $\T$ is a \emph{strongly continuous} representation, i.e. the function $x_\T: \G\to\X$ given by
$x_\T(g) =\T(g)x$ is continuous for every $x\in\X$, and $\nu(g) \geq \norm{\T(g)}$ for a weight $\nu$, then
a non-degenerate Banach $L_\nu(\G)$-module structure on $\X$ satisfying assumption \ref{assu} may be defined via
\begeq\label{scdef}
fx = \T(f)x = \int_\G  f(g)\T(-g)x dg,\quad f\in L_\nu(\G), \ x\in\X.
\eq
It is easy to show (\cite[Lemma 2.2]{BK05}) that every non-degenerate $L_\nu(\G)$-module has at most one representation
associated with it. Therefore, if $\T$ is strongly continuous, the module structure has to be defined via (\ref{scdef}). If $\T$ is not strongly continuous, we shall make use of the submodule $\X_{c}\subset \X$ of $\T$-\emph{continuous} vectors of $\X$. We let $x\in\X_c$ if the function
$x_\T$ is continuous.

The symbol $\T(f)$ in \eqref{scdef} may seem like an abuse of notation. We justify this by regarding $\T$ also as a representation of the algebra $M_\nu(\GG)$ of finite complex Borel measures on $\GG$ with respect to convolution and 
the norm 
\[\|\mu\|_\nu = \int_\GG \nu(g) d|\mu|(g).\]
The group $\GG$ is then identified with Dirac measures in $M_\nu(\GG)$ in an appropriate way  and $L_\nu(\GG)$ -- with the algebra of measures in  $M_\nu(\GG)$ that are absolutely continuous with respect to the Haar measure on $\GG$. For $x \in\X_c$ we then have
\begeq\label{tmes}
\mu x = \T(\mu)x = \int_\G  \T(-g)x d\mu(g),\quad \mu\in M_\nu(\G), \ x\in\X,
\eq
and $\|\mu x\|\le\|\mu\|_\nu\|x\|$. We refer to \cite{BK05} for more details.

\brem\label{tenp}
We observe that if $\G=\RR^d$ and we have an $L_1(\RR^d)$-module structure associated with a representation $\T:\RR^d\to B(\X)$ we also have $d$ module structures over $L_1(\RR)$ associated with the representations $\T^{(k)}: \RR\to B(\X)$, $k=1,\ldots, d$, given by
\[\T^{(k)}(t) x = \T(0,\ldots, 0,t,0,\ldots, 0)x,\ x\in\X,\]
where $t$ is the $k$-th component of the vector. Obviously, in this case,
\[\T(t_1,\ldots, t_d) = \prod_{k=1}^d \T^{(k)}(t_k)\quad\mbox{and}\quad
\T(f_1\otimes\ldots\otimes f_d) = \prod_{k=1}^d \T^{(k)}(f_k),\]
$t_k\in\RR$, $f_k\in L_1(\RR)$, $k=1,\ldots, d$. A similar observation is clearly valid in the case of tensor-product weights and Cartesian products of arbitrary LCA-groups.
\erem

\subsection{Spectral properties} The role of memory of operators on Banach modules is fulfilled by the notion of the Beurling spectrum. We begin by defining this notion for vectors in Banach modules. 

\bd
The \emph{Beurling spectrum} of a set $N \subseteq (\X,\T)$ is the subset $\L(N)= \L(N,\T)$ of the dual group $\Gg$ the compliment of which is given by
\[\{\g\in\Gg: \mbox{ there is } f\in L_\nu(\GG) \mbox{ such that } \hat f(\g) \ne 0 \]
\[\mbox{ and } fx =\T(f)x =  0 \mbox{ for all } x\in N\}.\]
\ed
When $N = \{x\}$ is a singleton, we shall write $\Lambda(x)$ instead of $\L(\{x\})$.
Given a closed set $\s\subset \Gg$ we shall denote by $\X(\s)$ the (closed) \emph{spectral submodule} of
all vectors $x\in\X$ such that $\L(x)\subseteq \s$. The symbol $\X_{Comp}$ will stand for the set
of all vectors such that $\L(x)$ is compact.

\brem\label{otherspec}
The notion of Beurling spectrum has a lot of very different guises (see \cite{B04, BR75, D56, EN06}  among many others). All of these were shown to be equivalent at the level of generality of this paper \cite{B04}; we especially note the equivalent notions of the spectrum used in \cite{GK10} and in
\cite[Theorem XI.11.24]{DS88II}. 
\erem

Observe that if $\X$ is an appropriate space of functions on $\GG$ and $\T=T$ is defined via
\eqref{trans1}, then $T(f)x = f*x$ and $\L(x) = \supp \hat x$, $x\in\X$, $f \in L_\nu(\GG)$, possibly, in the sense of distributions.
Alternatively, if $\X$ is an appropriate space of functions on $\Gg$ and $\T=M$ is defined via
\eqref{mod1}, then $M(f)x = \hat f x$ and $\L(x) = \supp x$, $x\in\X$, $f \in L_\nu(\GG)$.
If $\T = \T_\X$ is defined via a resolution of the identity as in \eqref{exres}, then
\[\L(x) = \{n\in\ZZ: \cP_n x\ne 0\}.\]

In the next lemma we present basic properties of the Beurling spectrum that will be heavily used throughout the paper. We refer to \cite{B04, BK05} and references therein for the proof.

\bl\label{sprop}
Let $\X$ be a non-degenerate Banach $L_\nu(\GG)$-module with the structure associated with a 
representation $\T$. 
Then 
\begin{description}
\item[ (i)]   $\L(M)$ is closed for every $M \subseteq \X$ and $\L(M)=\emptyset$ if and only if $M =\{0\}$;
\item[(ii)] $\L(Ax + By) \subseteq \L(x) \cup \L(y)$ for all $A$, $B \in B(\X)$ that commute with all 
operators $\T(f)$,
$f \in L_\nu(\GG)$; 
\item[(iii)] $\L(f x)\subseteq (\supp \hat f)\cap \L(x)$ for all $f \in L_\nu(\GG)$ and $x\in\X$;
\item[(iv)] $f x=0$ if $(\supp \hat f)\cap\L(x)=\emptyset$, where $f \in L_\nu(\GG)$ and $x\in\X$;
\item[(v)] $f x = x$ if $\L(x)$ is a compact set, and $\hat f \equiv 1$ in some neighborhood of $\L(x)$, $f\in L_\nu(\GG)$, $x\in\X$;
\item[(vi)] if $M_0$ is dense in $M\subseteq X$, then $\L(M) = \overline{\bigcup_{x\in M_0}\L(x)}$.
\end{description}
\el
From the above lemma and examples we see that the Beurling spectrum $\L(x,\T)$ is, indeed, well suited to represent the support of the ``Fourier transform'' of the function $x_\T(g)=\T(g)x$. The goal of this paper is to study the behavior of the decay of this ``Fourier transform'' (memory decay in the case of operators), relate it to
smoothness properties of the function $x_\T$, and, in particular, to prove Wiener-type theorems for $x$ when it is an invertible linear operator or an element in a Banach algebra. Appropriate module structures on the spaces of operators will be discussed in Section \ref{Clop}. Presently we shall introduce the notions of approximate identities and $\g$-nets which are essential for our analysis of the memory decay.

\subsection{\texorpdfstring{Approximate identities and $\g$-nets}{Approximate identities and gamma-nets}}
Bounded approximate identities (b.a.i.) and $\g$-nets are often used to obtain ergodic theorems in Banach modules \cite{B78, B88}. We shall use \bai to approximate $\T$-continuous vectors with vectors that have compact Beurling spectrum and $\g$-nets to approximate the ``$\g$-Fourier coefficient'' of $x_\T$. We index all our \bai and $\g$-nets by some net $\Omega$.

\bd
 A bounded net $(f_\a)$ in $L_\nu(\GG)$ is called a \emph{bounded approximate identity} (b.a.i) in the algebra $L_\nu(\GG)$ if the following two conditions hold:
 \begin{description}
\item[ (i)] $\hat f_\a(0) = 1$ for all $\a\in\Omega$. 
\item[(ii)] $\lim f_\a*f = f$ for all $f\in L_\nu(\GG)$.
\end{description}
\ed

The functions
\begeq\label{znet}
f_\a(t) = \frac1{\pi^d}\prod_{k=1}^d \frac{\a_k}{t_k^2+\a_k^2},
\eq
give an example of a \bai in $L_1(\RR^d)$ when $\a \to 0^+$. A construction of \bai in $L_1(\GG)$ with the compact support of the Fourier transform is presented, for example, in \cite{BK05}. Existence of \bai in $L_\nu(\G)$ is, to the best of our knowledge, an open question; we refer to \cite[Lemmas 3.3.10, 3.4.9]{B04} and \cite{L85, SW71} for more details on the subject. Here we will note only that the algebra $L_\nu(\RR^d)$ with a non-quazi-analytic weight $\nu$ contains a \bai if $\nu$ is majorized  by a monotonic non-quazi-analytic weight $\tilde\nu$, that is 
$\tilde\nu(u) \le \tilde\nu(v)$ whenever $|u|_2\le|v|_2$, $u,v\in\RR^d$.


The following useful result is well-known \cite[and references therein]{B04, BK05, HR79}.  There
by $\X_\Phi$ we mean the submodule of vectors admitting factorization: $\X_\Phi = \{x = fy\in\X, y\in\X, f\in L_\nu(\GG)\}$.

\bt\label{CoHe}
Let $(\X,\T)$ be a non-degenerate Banach $L_\nu(\GG)$-module. Then
$\X_{Comp} \subseteq \X_\Phi \subseteq {\X_{c}}$.
Moreover, if the algebra $L_\nu(\GG)$ contains a \bai then
\[\X_c = \overline{\X_{Comp}} = \{x = \lim_\a f_\a x,\ f_\a \mbox{ -- a \bai in } L_\nu(\GG)\}\]
and if there is a \bai with $\|f_\a\|\le 1$ then also $\X_c = \X_\Phi$.
\et

We will need the following class of subsets of the module $\X_\Phi$.

\bd\label{modorbit} 
Given $x\in (\X, \T)$ the \emph{modular orbit} $\Omega(x)\subset \X_\Phi$ is the set
$\Omega(x) =\{fx,\ f\in L_\nu(\GG)\}$.
\ed

Observe that by Domar's Theorem \cite{D56}, \cite[Ch. 5.1]{L85} we have that the set
$\{f\in L_\nu(\GG):\, \supp\hat f \mbox{ compact}\}$ is dense in $L_\nu(\GG)$ and, therefore,
due to Lemma \ref{sprop}, for any open neighborhood $\mathcal U$ of $\L(x)$ we have
\begeq\label{goodf}
\overline{\Omega(x)} = \overline{\{fx,\ f\in L_\nu(\GG), \mbox{ with }\supp\hat f\subset \mathcal{U}  \mbox{ compact}\}}. 
\eq

\bl\label{goodx}
If $x\in\X_c$ and $L_\nu(\GG)$ contains a \bai then $x\in \overline{\Omega(x)}$.
\el

\bpf
Follows from Domar's theorem and Theorem \ref{CoHe}.
\epf

\bd\label{gnet}
 For $\g\in\Gg$ a bounded net $(f_\a)$ in the algebra $L_\nu(\GG)$ is called a $\g$-net if the following two conditions hold:
\begin{description}
\item[ (i)] $\hat f_\a(\g) = 1$ for all $\a$;
\item[(ii)] $\lim f_\a * f = 0$ for every $f \in L_\nu(\GG)$ with $\hat f(\g) = 0$.
\end{description}
\ed

The functions \eqref{znet} provide an example of a $0$-net in $L_1(\RR^d)$ as $\a\to \infty$. Other
examples of $\g$-nets in $L_1(\GG)$ are furnished in \cite{BK05}. Unfortunately, existence of such nets in $L_\nu(\GG)$ presents a formidable restriction on the growth of the weight $\nu$. In particular, in case $\G=\RR$ it is known \cite{B04} that $\g$-nets do not exist already when $\nu$ is a linear weight. 

The following proposition is a version of the results stated in \cite{B78, BK05}.

\bp
Assume that $(f_\a)$ is a $\g$-net and $x\in\X$ satisfies $0\ne x_\g = \lim f_\a x$. Then
$\L(x_\g) = \{\g\}$ and the limit does not depend on the particular choice of the $\g$-net.
\ep

The above proposition shows why $\g$-nets are used to approximate ``Fourier coefficients'' of the function $x_\T$. Of course, in general the limit $x_\g$ does not exist. That is why in the following section (see \eqref{triangled}) we use elements
of a $\g$-net to define ``parts'' of a vector with the Beurling spectrum in a neighborhood of $\g$.

The following result uses  $\g$-nets in its proof and allows us to sharpen some of the properties of the Beurling spectra we mentioned in Lemma \ref{sprop}. 



\bl\label{ssprop}\cite[Lemma 3.7.32]{B04}.
Let 
$\X$ be a non-degenerate Banach $L_1(\GG)$-module with the structure associated with a 
representation $\T$. 
Then 
\begin{description}
\item[(i)] $f x=0$ if $(\supp \hat f)\cap\L(x)$ is countable and $\hat f(\g) = 0$ for all $\g \in (\supp \hat f)\cap\L(x)$, 
$f \in L_1(\GG)$, $x\in\X$;
\item[(ii)] $f x = x$ if $\L(x)$ is a compact set, the boundary of $\L(x)$ is countable, and $\hat f \equiv 1$ on  $\L(x)$, $f\in L_1(\GG)$, $x\in\X$.
\end{description}
\el


\subsection{Extended module structure}\label{closed}

In the following section we define various classes of spectral decay for vectors in Banach modules. To characterize some of these classes  we need to extend the Banach module structure on $\X$ by means of a closed operator functional calculus.

In the following definition we use the spaces
$\lnhat(\G)=\{\hat\varphi$: $\varphi\in L_\nu(\G) \}$ and
$\lnloc(\G)=\{h$: $\Gg \to \CC$ such that  
$h\hat\varphi\in\lnhat(\G)$  for all $\varphi\in L_\nu(\G)$
 with 
$\supp{\hat\varphi}$ compact$\}$. Observe that $\lnhat(\G)\subset\lnloc(\G)$.

For $h\in\lnloc(\G)$ we define a (closed) operator $\chT(h) = h\diamond : D(h)=D(h,\T)\subseteq\X\to\X$ in the following way.
First, let $x\in\X_{comp}$ and 
\begeq\label{cl}
\chT(h)x = h\diamond x : = (h\hat\varphi)^{\vee} x = \T((h\hat\varphi)^{\vee})x,
\eq
where $\varphi\in L_\nu(\G)$ is such that $\supp \hat\varphi$ is compact and  $\hat\varphi\equiv 1$ in a neighborhood of $\L(x)$.
The vector $\chT(h)x$ is well-defined in this way because it is independent of the choice of $\varphi$. Indeed, if $\phi\in L_\nu(\G)$ is another function with the same properties, then
$h\hat\phi-h\hat\varphi\in\lnhat(\G)$ is $0$ in a neighborhood of $\Lambda(x)$ and, hence,
\[ \T((h\hat\phi)^{\vee})x= \T((h\hat\varphi)^{\vee})x\]
by \eqref{cl} and Lemma \ref{sprop}. The same lemma ensures also that $h\diamond x\in\X_{comp}$ whenever $x\in\X_{comp}$.
Next, we extend the definition of $\chT(h)$ by taking the closure of the just defined operator on $\X_{comp}$. In other words,
if $x_n\in \X_{comp}$, $n\in\NN$, $x = \lim\limits_{n\to \infty} x_n$, and
$y = \lim\limits_{n\to \infty} h\diamond x_n$ exists, we let $h\diamond x = y$.

In the following lemma we show that the above definition makes sense.

\bl
The operator $\chT(h) = h\diamond : D(h)=D_{\T}(h)\subseteq\X\to\X$ 
is a well-defined closed operator.
\el

\bpf

We need to prove that the restriction of $\chT(h)$ onto $\X_{Comp}$ is closable. 
In particular, we need to show that if
$x_n\in \X_{comp}$, $n\in\NN$, $ \lim\limits_{n\to \infty} x_n = 0$, and
$y = \lim\limits_{n\to \infty} h\diamond x_n$ exists, then $y=0$. Consider an arbitrary 
$f\in L_\nu(\G)$ with $\supp\hat f$ compact and let $\varphi\in L_\nu(\G)$ satisfy
$\hat\varphi \equiv 1$ in a neighborhood of $\supp \hat f$. Let also $\varphi_n\in L_\nu(\G)$ be such that $\hat\varphi_n \equiv 1$ in a neighborhood of $\L(x_n)$, $n\in\NN$. Then
\begin{equation*}
\bs 
fy & = \lim_{n\to\infty} f(h\diamond x_n) = \lim_{n\to\infty} (h\hat f\hat\varphi_n)^\vee x_n
= \lim_{n\to\infty} (h\hat\varphi_n)^\vee (fx_n) 
\\ &= \lim_{n\to\infty} h\diamond (fx_n) 
=\lim_{n\to\infty} (h\hat \varphi)^\vee (fx_n) =0.
\end{split}
\end{equation*}
By Remark \ref{cnepr}, this implies $y = 0$.
\epf

From the definition we immediately have $\T(f) = \chT(\hat f)$, $f\in L_\nu(\G)$.
We also observe the following slightly less obvious properties of the extended module structure.

\bp\label{spro} If $g\in\G$, $f\in L_\nu(\G)$, $h\in \lnloc(\G)$, and $x\in D(h)$, then
$\T(g)x$, $\T(f)x \in D(h)$ and we have $\T(g)(h\diamond x) = h\diamond (\T(g)x))$ and
$\T(f)(h\diamond x) = h\diamond (\T(f)x))$.
\ep

\bpf
In case $x\in\X_{Comp}$, we get $\T(g)x$, $\T(f)x \in\X_{Comp}$ and the assertions are obvious. Otherwise, let $x =\lim\limits_{n\to\infty} x_n$ with $x_n\in\X_{Comp}$, $n\in\NN$. Then,
by the definition of $\chT(h)$ we have 
\[h\diamond (\T(g)x)) = \lim_{n\to\infty} h\diamond(\T(g)x_n) = \lim_{n\to\infty} \T(g)(h\diamond x_n)
=\T(g)(h\diamond x)\] 
and similar equations hold for $\T(f)(h\diamond x)$.
\epf

Some basic properties of the Beurling spectrum outlined in Lemmas \ref{sprop}, \ref{ssprop} remain valid for the extended module structure in the following way.

\bl\label{sssprop}
Assume $\L(x)$ is compact, has countable boundary, and $\varphi, \psi \in \lnloc(\G)$ satisfy 
$\varphi\equiv\psi$ on $\L(x)$. Then $\chT(\varphi)x = \chT(\psi)x$.
\el 

\bpf
Let $f\in L_1(\G)$ satisfy $\hat f\equiv 1$ in a neighborhood of $\L(x)$. Applying
Lemma \ref{ssprop} to the function $h =  (\hat f(\psi-\varphi))^\vee \in L_1(\G)$ we get
$0 = h x = \chT(\psi)x - \chT(\varphi)x$.
\epf

The following proposition will play an important role in estimating the spectral decay.

\bp\label{inmeas}
Assume that $h\in\lnloc(\G)$ is such that \[\hat h(\g) = \sum_{n\in\ZZ} c_n \g(g_n),\] $g_n\in\G$,
$n\in\ZZ$, and $\sum_{n\in\ZZ} |c_n|\nu(g_n) < \infty$. Then
\[\chT(h)x  = \sum_{n\in\ZZ} c_n \T(g_n)x = \T(\mu)x, \ x\in\X_c,\]
where the measure $\mu \in M_\nu(\GG)$ is given by $\mu = \sum_{n\in\ZZ} c_n \d_{-g_n}$, $\d_g$ -- Dirac measure concentrated at $g\in\GG$.
\ep

\bpf
Assume  $x\in \X_{Comp}$ and $f\in L_\nu(\G)$ is such that $fx = x$ and $\supp \hat f$ is compact. Then, using \eqref{redef},
\[\left(\chT(h) - \sum_{n\in\ZZ} c_n \T(g_n)\right)x =\left(\chT(h) - \sum_{n\in\ZZ} c_n \T(g_n)\right)\T(f)x=\]
\[\sum_{n\in\ZZ} c_n \chT(\hat f(\g)\g(g_n))x - \sum_{n\in\ZZ} c_n \T(T(g_n)f)x = 0,\]
and the proposition follows from the definition of $\chT(h)$ and Theorem \ref{CoHe}.
\epf

The following two types of functions in $\widehat{L_1}^{loc}(\RR^d)$ are especially useful. For $\l=(\l_1,\ldots,\l_d)\in\RR^d$ and $\a=(\a_1,\ldots,\a_d)\in {\RR^d}$ we let
\begeq\label{expfun}
\eet(\l) = e^{\a\cdot\l} = \exp\left(\sum\limits_{k=1}^d{\a_k  \l_k}\right)
\eq
\begeq\label{modexp}
\bh_\a (\l) =e^{\a\cdot|\l|} =  \exp\left(\sum\limits_{k=1}^d{\a_k | \l_k|}\right),
\eq
Observe that if  $\a\in\RR^d_+$, i.e. $\a_k>0$, $k=1,\ldots, d$, then $\bh_{-\a} = \hat f_\a$, where
$f_\a\in L_1(\RR^d)$ is the function introduced in \eqref{znet}, and for all $x\in D(\bh_\a)$ 
\begeq\label{fheq}
\bh_{-\a}\diamond(\bh_\a\diamond x) =f_{\a}(\bh_\a\diamond x) = \bh_\a\diamond(f_\a x) = x.
\eq
Observe also that if 
$\a=(0,\ldots,0,\a_k,0,\ldots,0)$, the following relation is valid:
\begeq\label{releh}
\eet+\textbf{e}_{-\a} = \bh_\a+\bh_{-\a},
\eq
and, hence, using $\bh_{-\a} \in \widehat{L_1}(\RR^d)$,
\begeq\label{dreleh1}
D(\eet)\cap D(\textbf{e}_{-\a})\subseteq D(\bh_\a),\ \a\in\RR^d_+.
\eq

Another family of functions we shall use is given by $e_z(\l) = e^{iz\cdot \l}$, $z\in \CC^d$. Observe that $e_{i\a} = \textbf{e}_{-\a}$, $\a\in\RR^d$, and, hence, Propositions \ref{spro} and \ref{inmeas} imply
\begeq\label{ee}
D(e_z) = D(\textbf{e}_{-\a}), \chT(e_z) = \chT(\textbf{e}_{-\a})\T(t),\ z = t+i\a,\ t,\a\in\RR^d.
\eq

\subsection{\texorpdfstring{Generators of $L_1(\RR)$-modules}{Generators of L1(\RR)-modules}}\label{gener} 
Consider $L_1(\RR)$-module $(\X,\T)$, satisfying $\X =\X_c$, that is the module structure is associated with a
 strongly continuous isometric representation $\T: \RR \to B(\X)$. In this section we consider the properties of the (infinitesimal) generator $iA$ of the one-parameter group $\T$ \cite{EN06}.
The operator $A$ is called the generator of the module $\X$.

\bt\label{spmap}
Assume $\L(\X)$ is compact and nonempty. Then the generator $A\in B(\X)$ and
we have the spectral mapping formula $\s(A) = \L(\X)$. Moreover, the function $\T:\RR\to B(\X)$ admits a holomorphic extension to an entire function $\T(z) = e^{izA} = \T(t)e^{-\a A} = e^{-\a A}\T(t) = \check\T(e_z)$, where $e_z(\l) = e^{iz\l}$, $z = t+i\a\in\CC$, and
\[\|\T(t+i\a)\| =\|\T(i\a)\| \le  \max_{\l\in\L(\X)}e^{-\a \l}, \ t,\a\in\RR.\]
\et

\bpf
Most of the assertions of this theorem are well known or nearly obvious (see, e.g.,  \cite{B79, B04}, \cite[Theorems 3.7, 3.8]{BK05} and references therein). In this paper, we shall prove only the norm estimate. We may assume that  $\a \in \RR_+$ and $a = |\inf \L(\X)| \ge |\sup \L(\X)|>0$ without loss of generality (otherwise we can consider the representation $\widetilde{\T}(t) = \T(-t)$). Let us estimate the norm
$\|\T(i\a)\|$.
Consider
the $4a$-periodic function $\eta$ defined on the interval $[-3a,a]$ via
\begeq\label{curvehat}
\eta(\l) = \left\{
\begin{array}{cl}
e^{-\a\l}, & \l\in[-a,a],\\
e^{\a (2a+\l)}, & \l\in[-3a,-a).
\end{array}
\right.
\eq
Computing the Fourier coefficients $c_n$, $n\in\ZZ$, of the periodic even function $\eta(\cdot - a)$ we get
\begeq\label{fcoef}
c_n = 2\int_0^{2a} e^{\a(a-\l)}\cos \frac{\pi n \l}{2a}d\l=\frac{8\a a^2 e^{\a a}}{4\a^2a^2 +\pi^2n^2}(1 - e^{-2\a a}\cos\pi n)>0.
\eq
Hence, in view of Proposition \ref{inmeas}, $\|\chT(\eta)\| \le 
e^{\a a}$.
Since $\eta \equiv e_{i\a}$
on  $[-a,a]\supseteq\L(A)$, Lemma \ref{sssprop} ensures $\|\T(i\a)\| =\|\chT(e_{i\a})\| = \|\chT(\eta)\|$ and the theorem is proved. 
\epf

\brem
It is not hard to show that the inequality in the norm estimate is, in fact, an equality. We do not need this fact in this paper and, therefore, will refrain from proving it.  
\erem

The following corollary is immediate. There we do not assume that $\L(\X)$ is compact.

\bc\label{Bineq}
Assume $x\in(\X,\T)$ and $\L(x)$ is compact and nonempty. Then the function $x_\T$: $\RR\to \X$ given by $x_\T(t) = \T(t)x$ admits a holomorphic extension to an entire function and 
\[\|x_\T(i\a)\| = \|e^{-\a A}x\| = \|\chT(e_{i\a})x\| \le \|x\|\max_{\l\in\L(\X)}e^{-\a \l},\ \a\in\RR.\]
\ec

\bpf
Lemma \ref{sprop} implies that $\L([x]) = \L(x)$, where $[x]\subseteq\X$ is the submodule generated by $x$, i.e. the smallest
submodule of $\X$ containing $x$. It remains to apply Theorem \ref{spmap} to $[x]$.
\epf

\brem
The above result leads to Bernstein-type inequalities that go back at least to \cite{B79} (see also \cite{BK05, GK10}). 
Note that if $x\in\X_{Comp}$ then $x\in D(A)$ and $\|Ax\| \le \rho(x)\|x\|$, where $\rho(x) = \max\limits_{\l\in\L(x)} |\l|$. Hence, if $\L(\X)$ is compact, then $A$ is a bounded operator and
$\|A\| = \rho(A)$, where $\rho(A)$ is the spectral radius of $A$, i.e. $\rho(A) = \max\{|\l|, \l\in\s(A)\}$.
\erem

\brem\label{multig}
In case $\GG=\RR^d$ we will consider \emph{multigenerators} of $L_1(\RR^d)$-modules (see Subsection \ref{sclass}). These are $d$-tuples comprised of the associated generators of $L_1(\RR)$-modules from the Remark \ref{tenp}. We also note, in passing, that for the purposes of this paper (proving Wiener-type results) it is not necessary to consider representations $\T$ that are not strongly continuous (in which case the generator may not be densely defined). In view of Theorem \ref{CoHe}, all one would need to do if $\T$ is not strongly continuous to begin with, is to restrict attention to the submodule $\X_c$. As we shall se below, the results on inverse closedness will then be obtained without any loss of generality (compare with \cite{GK10}). 
\erem

\section{Different classes of spectral decay}\label{Cla}
In this section we introduce classes of elements with different spectral decay. 
Unless stated otherwise, in this section we assume that $\G=\RR^d$ and $\nu\equiv 1$.

\subsection{Exponential decay, one-sided exponential decay and causal classes}\label{exclass}
A major role in this paper is played by classes of elements with exponential spectral decay.

\bd\label{expdecay}
We say that $x\in(\X,\T)$ has \emph{exponential spectral decay of type} $\a\in\RR^d_+$ if $x\in D(\bh_\a)$. Similarly, $x$ has \emph{one-sided exponential spectral decay of type} $\a\in\RR^d$ if $x\in D(\eet)$.
\ed

In the following lemmas we provide useful characterizations of the newly defined classes. In the proof we use the classes of causal (and anticausal) elements which we define as in \cite{BK05}.

\bd
The submodules $\C_+(\X)$ and $\C_-(\X)$ of \emph{causal} and \emph{anticausal} elements in $\X$, respectively, are defined via
\[\C_\pm(\X) = \{x\in\X:\, \L(x)\subseteq \overline{\RR^d_\pm}\}.\]
\ed

We begin the study of the classes with a characterization of $D(\eet)$ in case $d=1$. In what follows we
consider a family of functions
$\phi_N \in L_1(\RR^d)$, $N\in\NN$, such that $\norm{\phi_N}=1$, 
$\hat\phi_N$ is continuous, $\supp\hat\phi_N\subset [-N,N]^d$, and $\phi_{N,n}(s) = \phi_N(s)e^{ i Nn\cdot s}$ so that $\hat\phi_{N,n}(\l) = \hat\phi_N(\l - Nn)$, $n\in\ZZ^d$.
In this paper we routinely use specific
 functions $\phi_N$, $N\in\NN$, defined via
\begeq\label{triangle}
\hat\phi_N(\l) = \hat\phi(\frac\l N),\quad
\hat\phi(\l) = \hat\phi_1(\l)=(1-|\l|)\chi_{[-1,1]}(\l),
\eq
where $\chi_S$ is, as usually, a characteristic function of the set $S$. It should be understood, however, that other functions with the specified properties may be used just as well.

\brem\label{Schwrep}
We can replace the function $\phi$ in \eqref{triangle} with a function $\varphi$ that is in the Schwartz class $\mathcal S$. As long as we have $\supp \hat\varphi\subseteq [-1,1]$ and
\[\sum_{n\in\ZZ} \hat \varphi (\l - n) \equiv 1,\]
all of the results in this paper remain valid.
\erem

\bl\label{1side}
Let $x\in\X_c$, $N\in\NN$, and $\phi_{N,n}$,  $n\in\ZZ$, be defined by \eqref{triangle}. 
Then the following are equivalent:
\begin{description}
\item[(i)] $x\in \bigcup_{\a > 0} D(\emph{\textbf{e}}_{\pm\a})$;
\item[(ii)] there are $M> 0$ and $\g\in [0,1)$ such that $\|\phi_{N, \pm n} x\|\le M\g^{Nn}$ for  all $n \in \NN$ (equivalently, for all $n\in\ZZ$ with $n\ge k$, $k\in\ZZ$);
\item[(iii)] the function $x_\T: \RR\to\X$ given by
$x_\T(t) =\T(t)x$, $t\in\RR$, admits a holomorphic extension to the strip
\[\CC_{\mp\d} = \{z\in\CC: 0 < \mp {\Im}m z < \d\}\]
for some $\d \in (0,\infty]$, which is continuous and uniformly bounded in $\overline{\CC_{\mp\d}}$.
\end{description}
\el

\bpf
(i)$\implies$(ii). Let us assume that $x\in  D(\eet)$ for some $\a > 0$ and proceed to estimate $\|\phi_{1,n}x\|$ for $n\in\NN$. Using Lemma \ref{ssprop} and Proposition \ref{spro}, we get 
\[\|\phi_{N,n}x\| =  \|f_\a(\eet\diamond(\phi_{N,n}x))\| = \|(f_{\a}*\phi_{N,n})(\eet\diamond x)\| \le \|f_{N,n}\|\|\eet\diamond x\|,\]
where $f_{N,n}\in L_1(\RR^d)$ with $\hat f_{N,n} (\l)= \frac{e^{\alpha N(1-n)}}{\ee(N(\l-n+1))}\in \widehat{L_1}(\RR^d)$, $n\in\NN$. Since $f_{N,n}(s)e^{ i N(n-1)s}$ is positive for all $x\in\RR$, $\|f_{N,n}\| = \hat f_n(n-1) = e^{\alpha N(1-n)}$ and we deduce that 
$\|\phi_{N,n}x\| \le  e^{\alpha N(1-n)}\|\eet\diamond x\|$ for $n \in \NN$. Similar argument works for
$x \in D(\textbf{e}_{-\a})$ and we obtain
\begeq\label{eexpest}
\|\phi_{N,\pm n}x\| \le  e^{\alpha N(1-n)}\|\textbf{e}_{\pm\a}\diamond x\|,\ n\in \NN.
\eq

(ii)$\implies$(iii). Assume now that  $\|\phi_{N, -n} x\|\le M\g^{Nn}$, $n\ge -1$ (the other case is handled similarly). Then, we have $x = x_+ + x_-$, where 
$x_- = \sum\limits_{n\ge -1} \phi_{N,-n}x$ with the series converging absolutely and $x_+= x - x_-$.  Let us prove that $x_+ \in \C_+$. Indeed,  let $\l\in [-Nk, -N(k-1))$, $k\in\NN$. Consider a function $f\in L_1(\RR)$ with $\hat f(\l)\neq 0$ and $\supp\hat f\subset  (-N(k+1), -N(k-1))$. Then $$fx_+ = fx - f\left( \sum\limits_{n\ge -1} \phi_{N,-n}x\right) = fx - \left(f * \sum\limits_{n= k-2}^{k+2} \phi_{N,-n}\right)x = 0.$$ Therefore, $\l\notin\L(x_+)$. By \cite[Lemma 8.2]{BK05} the function $(x_+)_\T$: $t\mapsto \T(t)x_+$ admits the desired extension to the halfplane $\CC_{\infty}$. Next, we observe that each vector $\phi_{N,n}x$, $n\in\ZZ$, has compact Beurling spectrum and, hence, Corollary \ref{Bineq} applies. In this way we obtain holomorphic extensions $(\phi_{N,-n}x)_\T$, $n\ge -1$, that satisfy
\[\|(\phi_{N,-n}x)_\T(t+i\a)\| \le e^{|\a| N(1-n)}\|\phi_{N,-n}x\| \le Me^{|\a| N} (\g e^{-|\a|})^{Nn},\ t,\a\in\RR.\]
Hence, the series 
\[(x_-)_\T(z) = \sum_{n\ge -1} (\phi_{N,-n}x)_\T(z)\]
converges absolutely and uniformly in the strip $\CC_{\d}$, where $0 < \d < -\ln \g$. Therefore,
$x_\T (z) = (x_-)_\T(z) +(x_+)_\T(z)$ is a well-defined bounded function which is holomorphic in $\CC_\d$ and
continuous in $\overline{\CC_\d}$.

(iii)$\implies$(i). Assume that $x_\T$ is well defined, holomorphic in $\CC_{-\d}$ for some $\d > 0$ as well as continuous and bounded in $\overline{\CC_{-\d}}$. Let $f_k$ be a b.a.i. in $L_1(\RR)$ with $\supp \hat f_k$ compact, $k\in\NN$. 
The functions $z\mapsto e_z\diamond (f_k x)$ and $z\mapsto f_k(x_\T(z))$ are both holomorphic  in $\CC_{-\d}$ and continuous and bounded in $\overline{\CC_{-\d}}$. Since they coincide on $\RR$, they coincide in $\CC_{-\d}$ as well. Hence, the limit 
\[\lim_{k\to\infty} e_z\diamond (f_k x) = \lim_{k\to\infty} f_k (x_\T(z)) = x_\T(z) \]
exists and, therefore, $x\in D(e_z) = D(\textbf{e}_{-\a})$ for all $z = t+i\a \in \CC_{-\d}$ by \eqref{ee}.  
\epf

Observe that the estimates in the first part of the above proof remain valid when $\eet$ is replaced
by $\ee$, that is
\begeq\label{expest}
\|\phi_{N,n}x\| \le  e^{\alpha N(1-|n|)}\|\ee\diamond x\|,\ |n|\ge1.
\eq
In addition, since
$\phi_N\equiv \phi_{N,0}\in L_1(\RR)$ and $\|\phi_N\| = 1$, we have that
\[\|\phi_Nx\|=\|\phi_{N,0}x\|\le\|x\|\le\|\bh_\a\diamond x\|, N\in\NN, \a>0,\]
due to \eqref{fheq}.

The above estimates indicate that the terminology in Definition \ref{expdecay} is appropriate. Moreover, together with \eqref{dreleh1} and Remark \ref{tenp} they imply
\begeq\label{dreleh}
D(\eet)\cap D(\textbf{e}_{-\a})= D(\bh_\a),\ \a\in\RR^d_+.
\eq

The following lemma is almost immediate now.

\bl\label{2side}
Let $x\in\X_c$, $N\in\NN$, and $\phi_{N,n}$,  $n\in\ZZ$, be defined by \eqref{triangle}. 
Then the following are equivalent:
\begin{description}
\item[(i)] $x\in \bigcup_{\a > 0} D(\emph{\textbf{h}}_{\a})$;
\item[(ii)] there are $M\ge 0$ and $\g\in [0,1)$ such that $\|\phi_{N, n} x\|\le M\g^{N|n|}$ for  all $n \in \ZZ$;
\item[(iii)] the function $x_\T: \RR\to\X$ given by
$x_\T(t) =\T(t)x$, $t\in\RR$, admits a holomorphic extension to the strip
\[\CC_{-\d,\d} = \{z\in\CC: -\d <  {\Im}m z < \d\}\]
for some $\d \in (0,\infty]$, which is continuous and uniformly bounded in $\overline{\CC_{-\d,\d}}$.
\end{description}
\el

\bpf
Assume $x\in D(\ee)$. In view of \eqref{dreleh}, Lemma \ref{1side} provides us with a holomorphic extension of $x_\T$ to $\CC_{\pm\d}$ for some $\d >0$.  The only thing that remains to prove is that this extension is holomorphic in a neighborhood of $\RR\subset \CC$. Observe that the family of functions $x_\a \equiv x_\T(\cdot + i\a): \RR\to\X$ is uniformly bounded and equicontinuous for $\a\in[-\d,\d]$. Hence, the extension $x_\T$ is holomorphic in $\CC_{-\d,\d}$ as a uniform limit of a sequence of entire functions $f_n x_\T$, where $(f_n)$ is a \bai in $L_1(\RR)$ with $\supp \hat f_n$ compact.
\epf

Remark \ref{tenp} allows us to extend Lemmas \ref{1side} and \ref{2side} to the multidimensional case.
Using the functions from \eqref{triangle}, we define $\phi_{N,a}^d\in L_1(\RR^d)$, $a\in\RR^d$, $N\in\NN$, via 
\begeq\label{triangled}
\hat\phi_{N,a}^d (\l) = \hat\phi_{N}^d(\l-Na), \ \mbox{ where }
\hat \phi_{N}^d (\l_1,\ldots,\l_d)= \prod_{k=1}^d \hat\phi_{N}(\l_k).  
\eq

We obtain the following characterizations of the classes with exponential spectral decay.

\bp\label{1expestprop}
Let $x\in\X_c$, 
$N\in\NN$, and $\phi_{N,n}^d$,  $n\in\ZZ^d$, be defined by \eqref{triangle} and \eqref{triangled}. Then the following are equivalent:
\begin{description}
\item[(i)]  $x\in \bigcup_{\a \in\RR^d_\pm} D(\emph{\textbf{e}}_{\a})$;
\item[(ii)]  there are $M\ge 0$ and $\g\in [0,1)$ such that $\|\phi_{N,\pm n}^d x\|\le M\g^{N|n|}$, $n\in\ZZ^d_+$; 
\item[(iii)] the function $x_\T: \RR\to\X$ given by
$x_\T(t) =\T(t)x$, $t\in\RR$, admits a holomorphic extension to the strip
\[\CC^d_{\mp\d} = \{z = (z_1, z_2, \ldots, z_d)\in\CC^d: 0 < \mp {\Im}m z_k < \d_k, k=1,2,\ldots, d\}\]
for some $\d = (\d_1, \d_2,\ldots, \d_d) \in \CC_+^d$, which is continuous and uniformly bounded in $\overline{\CC^d_{\mp\d}}$.
\end{description}
\ep

\bp\label{expestprop}
Let $x\in\X_c$, 
$N\in\NN$, and $\phi_{N,n}^d$,  $n\in\ZZ^d$, be defined by \eqref{triangle} and \eqref{triangled}. 
Then the following are equivalent:
\begin{description}
\item[(i)]  $x\in \bigcup_{\a \in\RR^d_+} D(\emph{\textbf{h}}_{\a})$;
\item[(ii)]  there are $M\ge 0$ and $\g\in [0,1)$ such that $\|\phi_{N,\pm n}^d x\|\le M\g^{N|n|}$;
\item[(iii)] the function $x_\T: \RR\to\X$ given by
$x_\T(t) =\T(t)x$, $t\in\RR$, admits a holomorphic extension to the strip
\[\CC^d_{-\d,\d} = \{z = (z_1, z_2, \ldots, z_d)\in\CC^d: -\d_k <  {\Im}m z_k < \d_k, k=1,2,\ldots, d\}\]
for some $\d = (\d_1, \d_2,\ldots, \d_d) \in \CC_+^d$, which is continuous and uniformly bounded in $\overline{\CC_{\mp\d}}$.
\end{description}
\ep


We will often use the following more precise estimates  that were obtained in the course of the proof: 
\begeq\label{expestd}
\|\phi_{N,n}^dx\|\le  e^{NM_n(\a)}\|\ee\diamond x\|, \ x\in D(\ee), \ n\in\ZZ^d,
\eq
\begeq\label{eexpestd}
\|\phi_{N,\pm n}^dx\|\le  e^{NM_n(\a)}\|\textbf{e}_{\pm\a}\diamond x\|, \ x\in D(\textbf{e}_{\pm\a}),\
n\in\ZZ^d_+,
\eq
where 
\begeq\label{mn}
M_n (\a)= \sum_{k:\ n_k\ne 0} \a_k(1-|n_k|), \ n=(n_1,\ldots,n_d)\in\ZZ^d.
\eq

Observe also that for any
$\a\in\RR^d_+$ we have $\X_{Comp} \subseteq D(\eet)$ and 
$\X_c\cap D({\textbf{e}}_{\beta}) \subseteq X_c\cap D({\textbf{e}}_{\a})$ whenever $\beta\in\RR^d_\pm$ and $\a-\beta \in\RR^d_\pm$. We also mention explicitly that extending \cite[Lemma 8.2]{BK05} via Remark \ref{tenp}, we see that for every
$\T$-continuous vector $x\in\C_+(\X)$ the function $x_\T$ has a bounded holomorphic extension to the halfspace $\CC_{\infty}^d$ given by 
\begeq\label{holex}
 x_\T(z) = \bh_{\a}\diamond (\T(t)x)= \textbf{e}_{-\a}\diamond (\T(t)x) = e^{iz\cdot\l}\diamond x,\ z = t+i\a\in\overline{\CC^d_{\infty}},
\eq
and a similar formula holds for $x\in \C_-(\X)\cap\X_c$.

\subsection{Wiener Class}\label{wiencl}

The analog of the classical Wiener class  $\W\subseteq \X$ in our setting is defined as follows: 
\begeq
\W =\W(\X) = \W(\X,\T):= \left\{x\in\X:\ \|x\|_\W =\int_{\RR^d}\|\phi^d_{1,a} x\|da< \infty \right\},
\eq
where the functions $\phi^d_{1,a}$ are given by \eqref{triangled}. Since the representation $M$ defined by \eqref{mod1} on $L_1(\RR^d)$ is strongly continuous,  the function $a\mapsto \phi^d_{1,a} x: \RR^d\to\X$ is continuous and, therefore, the finiteness of the integral is a reasonable way to define the corresponding type of spectral decay. Presently, we shall see that
$\W$ is a Banach space; for our purposes, however, it is more convenient to use equivalent
series norms rather than the integral norm in the definition. We define these norms via
\begeq\label{wiensumel}
\|x\|_{1,N} = 5^d \sum_{n\in\ZZ^d} \|\phi^d_{N,n} x\|,
\ N\in\NN.
\eq
The seemingly strange choice of the constant $5^d$ will become natural when we study the Wiener Class of operators.

\bp\label{preq}
The set $\W$ is a Banach space with respect to any of the equivalent norms in \eqref{wiensumel}.
Moreover,
\begeq\label{normeq}
\frac1{(2N+1)^d} \|x\|_{1,1} \le \|x\|_{1,N} \le 3^d\|x\|_{1,1},\ \mbox{and}
\eq
\begeq\label{normeq1}
 \|x\|_{\W} \le \|x\|_{1,1} \le 20^d\|x\|_{\W},\ x\in\W.
\eq
\ep
\bpf
We begin by showing  that $x\in\W$ if and only if $\|x\|_{1,1}<\infty$, obtaining \eqref{normeq1} in the process. The crucial observation  here is (for $d=1$)
that \[\phi_{1,a}x =\phi_{1,a}((\phi_{1,n-1}+\phi_{1,n}+\phi_{1,n+1}+\phi_{1,n+2})x),\] where
$n = \lfloor a \rfloor$ is the largest integer less than or equal to $a\in\RR$. A similar equality for $a\in\RR^d$ immediately implies
\[ \int_{\RR^d}\|\phi_{1,a}^d x\|da \le 4^d \sum_{n\in \ZZ^d} \|\phi_{1,n}^dx\| \le \|x\|_{1,1}.\]

Similarly, since
\[\|\phi_{1,n}^dx\| = \int_{n+[0,1]^d} \left\|\phi_{1,n}^d\left(\left(\sum_{k\in S_n}\phi_{1,a-k}^d\right)x\right)\right\|da,\]
where each $S_n$, $n\in\ZZ^d$, is the set of cardinality at most $4^d$ such that the Fourier transform of the sum is equal to 1 on $\supp \hat\phi_{1,n}^d$, we have
\[\|x\|_{1,1} \le 20^d  \int_{\RR^d}\|\phi^d_{1,a} x\|da.\]

Next, we prove the inequalities \eqref{normeq}. We use Lemma \ref{ssprop} and Remark \ref{tenp} in the usual way to obtain them.
The second of the estimates follows from
\begin{equation*}
\bs
\sum_{n\in\ZZ^d} \|\phi^d_{N,n}x\| &= \sum_{n\in\ZZ^d} \left\|\phi^d_{N,n}\left(\sum_{k\in\ZZ^d} \phi^d_{1,k}x\right)\right\| \\
&\le\sum_{k\in\ZZ^d} \sum_{n\in S_k}\|\phi^d_{N,n}\| \left\|\phi^d_{1,k}x\right\| \le 3^d \sum_{k\in\ZZ^d} \left\|\phi^d_{1,k}x\right\|,
\end{split}
\end{equation*}
where $S_k\subset \ZZ^d$ is a set of $d$-tuples  $n\in\ZZ^d$ such that the intersection of the supports
 of $\hat\phi^d_{N,n}$ and $\hat\phi^d_{1,k}$ has non-empty interior (all these sets have cardinality at most $3^d$).
 Similarly, the first inequality in \eqref{normeq} follows from
\begin{equation*}
\bs
\sum_{n\in\ZZ^d} \|\phi^d_{1,n}x\| &= \sum_{n\in\ZZ^d} \left\|\phi^d_{1,n}\left(\sum_{k\in\ZZ^d} \phi^d_{N,k}x\right)\right\| \\
&\le \sum_{k\in\ZZ^d} \sum_{n\in S^k}\|\phi^d_{1,n}\| \left\|\phi^d_{N,k}A\right\| = (2N + 1)^d\sum_{k\in\ZZ^d} \left\|\phi^d_{N,k}x\right\|,
\end{split}
\end{equation*}
where $S^k\subset \ZZ^d$ is a set of numbers  $n\in\ZZ^d$ such that the intersection of the supports
 of $\hat\phi^d_{1,n}$ and $\hat\phi^d_{N,k}$ has non-empty interior (all these sets have cardinality 
 $(2N+1)^d$).

The rest of the assertions in the proposition are now trivial.\epf

\brem\label{wc}
Observe that in view of Theorem \ref{CoHe} the definition of Wiener class implies $\W(\X)\subseteq \X_c$. This is why we decided to limit many of the auxiliary statements appearing in this 
paper to $\X_c$ even though in most cases they can be stated more generally. 
\erem

Let us present a few simple examples of the class $\W$.

\bex\label{ex1}
Let $\X = L^p(\RR^d)$, $1\le p \le\infty$, be the space of equivalence classes of $p$-integrable complex-valued functions, and $T$ and $M$ be the translation and modulation representations defined as in \eqref{trans1} and \eqref{mod1}. Let also $U$ be a representation based on a resolution of the identity $\mathcal{P} = \{P_k\}_{k\in\ZZ^d}$, as in \eqref{exres}. Then
\[
\bs
\W(\X, M) & = \left\{f\in\X:\ \sum_{n\in\ZZ^d} \left(\int_{\RR^d}\left|\hat\phi^d_{1,n}(t)x(t)\right|^pdt\right)^{\frac1p}<\infty\right\} \\
& =  \left\{f\in\X:\ \sum_{n\in\ZZ^d} \left(\int_{[-1,1]^d}\left|x(t-n)\right|^pdt\right)^{\frac1p}<\infty\right\}, 
\end{split}\]
 $1\le p< \infty$, with the obvious change for
$p=\infty$, 
is the Wiener amalgam space $L^{p,1}(\RR^d)$ \cite{FS85, H75, W32}. For $p=2$ we also have
$\W(L^2(\RR^d), T) = \widehat L^{2,1}(\RR^d) = \{x\in L^2(\RR^d):\ \hat x\in L^{2,1}(\RR^d)\}$.
In case of the representation $U$, we get $\W(\X, U) = \{x\in\X: \sum_{k\in\ZZ^d} \|P_kx\|<\infty\}$. In particular, if $(P_k x)(t) = \chi_{k+[0,1]^d}(t)x(t)$, we get  $\W(\X, U) = \W(\X, M) = L^{p,1}(\RR^d)$.
Alternatively, if $(P_k x)(t) = (\chi_{k+[0,1]^d}\hat x)^\vee(t)$, $x\in L^2(\RR^d)$, we get
 $\W(L^2(\RR^d), U) = \W(L^2(\RR^d), T)$. This shows that it may not be advantageous to consider 
 representations other than $U$ when $\X = L^p$. As a consequence, when we get to presenting 
 examples in Section \ref{examp} we concentrate primarily on Banach spaces other than $L^p$. Indeed,
 in case of more general Banach spaces, desired resolutions of the identity may not exist and the use of different representations becomes extremely beneficial.
 \eex

\subsection{Beurling Class} This class of elements is designed to encompass the classical
Beurling algebra of functions
\[A^*(\TT) = \left\{\sum_{n=-\infty}^\infty a(n)e^{in\xi}:\ \sum_{k=0}^\infty \sup_{|n|\ge k}|a(n)|<\infty\right\}\]
and was inspired by the work of Q.~Sun in \cite{S10ca}, see also \cite{BLT97}.

Again, we use the functions $\phi^d_{1,n}$, $n\in\ZZ^d$, defined by \eqref{triangled} and let
$\B = \{x\in \X:\ \|x\|_{\B}<\infty\}$, where
\begeq\label{bnorm}
\|x\|_{\B}= 
\sum_{k\in\ZZ^d}^\infty \max_{|n|_\infty\ge |k|_\infty} \|\phi^d_{1,n}x\|,
\eq
where $|k|_\infty = \max\{k_1,k_2,\ldots, k_d\}$, $k = (k_1,k_2,\ldots, k_d)$.

It is easily seen that \eqref{bnorm} defines a Banach space norm on $\B$ and for any $\a > 0$
we have
\begeq
D(\ee) \subseteq \B\subseteq \W.
\eq

\subsection{Sobolev-type classes}\label{sclass}
It is natural to consider weighted extensions of the classes introduced in this section. In general, we shall pursue such extensions in a sequel to this paper. Here we limit ourselves to Sobolev-type classes defined as follows (see Remark \ref{winsob}).

As usually, we consider a Banach space $\X$ equipped with a non-degenerate $L_1(\RR^d)$-module structure associated with a strongly continuous isometric representation $\T:\RR^d\to B(\X)$. We let $A = (A_1,\ldots, A_d)$ be the \emph{multigenerator} of the module $(\X,\T)$, i.e. each operator $iA_k$ is the infinitesimal generator of the one-parameter group $\T^{(k)}:\RR\to B(\X)$ from Remark \ref{tenp}, $k = 1,\ldots, d$. In the following definitions we use the standard multi-index notation, i.e. $\a = (\a_1,\ldots \a_d)\in \ZZ^d_+$,
$|\a| =|\a|_1 = \a_1+\ldots +\a_d$, and $A^\a = A_1^{\a_1}\ldots A_d^{\a_d}$.

\bd
For each $m\in\NN$ the \emph{Sobolev-Wiener class} $\W^{(m)}$ is given by
\[\W^{(m)} = \{x\in\X:\, x\in D(A^\a) \mbox{ and } A^\a x\in\W \mbox{ for all } |\a|\le m, \a\in\ZZ^d_+\}.\]
\ed

\bd
For each $m\in\NN$ the \emph{Sobolev-Beurling class} $\B^{(m)}$ is given by
\[\B^{(m)} = \{x\in\X:\, x\in D(A^\a) \mbox{ and } A^\a x\in\B \mbox{ for all } |\a|\le m, \a\in\ZZ^d_+\}.\]
\ed

\brem\label{winsob}
It can be easily shown that for $\a = (m,\ldots, m)\in \NN^d$ we have $D(A^{\a}) = D_\T(\varphi_{m})$, where
$\varphi_{m}(\l) = (\l_1\ldots \l_d)^{m}$, $\l = (\l_1,\ldots, \l_d)\in \RR^d$. Hence, $A\in \W^{(m)}$  if and only if
\[\int_{\RR^d}\|\phi^d_{1,a}A\|(1+|a|)^{m}da < \infty.\]
We refer to \cite{GK10, Kl11} for similar and other classes utilizing $D(A^\a)$.
\erem

\brem
It is not hard to see that in the case of Hilbert spaces or when $\X = C_{ub}(\RR^d,\HH)$ is a space of bounded uniformly continuous functions with values in a Hilbert space $\HH$, smoothness of the function $x_\T$ (the property $x\in D(A^\a)$) is closely related to the spectral decay of $x$ (see, e.g., \cite{GK10}). In particular, for $x\in C_{ub}(\RR,\HH)=\X$ we have that $x\in D(A)$ implies $x\in \W(\X, T)$, where $T$ is as in \eqref{trans1}. In general Banach spaces, however, this is not the case. Consider the following example. Let $\X = C_{ub}(\RR, C_{ub}(\RR))$ and $\{\a_n\}_{n\in\ZZ}$ be a null sequence ($\lim a_n = 0$) such that the series
$\sum\limits_{n\in\ZZ} \frac{\a_n}n$ diverges. Define $\varphi_n = \a_n\hat\phi_{\frac12, n}$ where $\phi_{\frac12,n}$ is as in \eqref{triangled} and let $x\in\X$ be the periodic function given by
\[x(t) = \sum_{n\in\ZZ} \frac1n\varphi_ne^{int},\  \varphi_n\in C_{ub}(\RR),\]
The function is well defined because the series converges unconditionally in $\X$. Moreover, the same is true for
\[x^\prime(t) = \sum_{n\in\ZZ} i\varphi_ne^{int},\]
so that $x\in D(A)$, $A =-i\frac{d}{dt}$. On the other hand, it is clear that $x\notin \W(\X,T)$. Hence, in general, the relation between smoothness of $x_\T$ and the spectral decay of $x$ is subtler.
Let us prove that $x\in D(A^2)$ implies $x\in \W(\X,\T)$ for general $\X$ and $\T: \RR\to B(\X)$. Let
$y\in D(A^2)$. Then there exists $x\in\X$ such that $y = (A^2+I)^{-1}x = fx$, where $\hat f(\l) = (\l^2+1)^{-1}$.  Hence, $\phi_{1,n}y = (f*\phi_{1,n})x$ satisfies $\|\phi_{1,n}y\| \le \|f*\phi_{1,n}\|\|x\| \le \frac{c}{n^2}$ for some $c\in\RR$ and, therefore, $y\in\W(\X,\T)$. Clearly, a small modification of the above argument shows
$D(A^{1+\varepsilon})\subseteq \W(\X,\T)$ for any $\varepsilon > 0$.
\erem

\brem\label{periodize}
We defined various classes of vectors with spectral decay assuming $\G = \RR^d$. Such a restriction, however, is less stringent than one may think. Indeed, if $\G = \TT^{d_1}\times\RR^{d_2}$ one can pass to a representation
on $\RR^{d_1+d_2}$ which would be periodic in the first $d_1$ variables and obtain the same kind of results  just as well. 
\erem

\section{Module structures and memory of linear operators}\label{Clop}

In view of our interest in Wiener-type results, we consider the classes defined in the previous section for Banach modules $L(\X, \Y)$ and $L(\Y, \X)$ of bounded linear operators between
complex Banach spaces $\X$ and $\Y$. To avoid confusion we shall use 
$\T_{\X\Y}: \GG\to B(L(\X,\Y))$ to denote the representation defining the module structure on
$L(\X,\Y)$ and $\T_{\Y\X}: \GG\to B(L(\Y,\X))$ -- for the representation on $L(\Y,\X)$. Since we are interested only in operators with some kind of memory decay we will only deal with operators in
$(L(\X,\Y))_c$ and $(L(\Y,\X))_c$. This allows us to use formula \eqref{scdef} and avoid tedious complications that we had to deal with in \cite{BK05}. 
We shall also always assume that the 
representations $\T_{\X\Y}$ and $\T_{\Y\X}$ are \emph{coupled}, i.e. 
the operators $A\in L(\X,\Y)$, $B\in L(\Y,\X)$ are inverses of each other if and only if
$\T_{\X\Y}(g)A$ and $\T_{\Y\X}(g)B$ have the same property for every $g\in\GG$.
Moreover, given three (not necessarily distinct) complex Banach spaces $\X$, $\Y$, $\Z$ and
Banach modules $L(\X,\Y)$, $L(\Y,\Z)$, and $L(\X,\Z)$ with structures associated with $\T_{\X\Y}$,
$\T_{\Y\Z}$ and $\T_{\X\Z}$, respectively, we assume that
\begeq\label{couple}
\T_{\X\Z}(g)(BA) = (\T_{\Y\Z}(g)B)(\T_{\X\Y}(g)A),\ g\in\GG,
\eq
for all $A\in L(\X,\Y)$, $B\in L(\Y,\Z)$. In particular, if $\X = \Y$, we assume that $\T_{\X\X}(g)$, $g\in\G$, are Banach algebra automorphisms of $B(\X)$.

If the spaces $\X$, $\Y$, and $\Z$ are themselves Banach modules, we shall use $\T_\X$, $\T_\Y$, and $\T_\Z$ to denote the corresponding representations.
However, in what follows,  we omit the indices of the representations if the choice is unambiguous.

\bex\label{typop}
The most typical example of module structures on $L(\X, \Y)$ and $L(\Y,\X)$ satisfying the above assumptions arises when $\X$ and $\Y$ are equipped with module structures associated with $\T_\X$ and $\T_\Y$, respectively. In this case, we let $\T_{\X\Y}(g)A = \T_\Y(g)A\T_\X(-g)$ and
$\T_{\Y\X}(g)B = \T_\X(g)B\T_\Y(-g)$, $g\in\GG$, $A\in L(\X,\Y)$, $B\in L(\Y,\X)$. Clearly, these are coupled representations and \eqref{couple} is satisfied in case we have three spaces and define the representations in this way. We also observe that when the representations $\T_\X$ and $\T_\Y$ are defined via resolutions of the identity as in \eqref{exres} and $A_{mn}$ are operator blocks
of an operator $A\in L(\X,\Y)$ as described in the introduction, then $k\notin \L(A,\T_{\X\Y})$ if and only if $A_{mn} = 0$ for all $m,n \in \ZZ^d$ such that $m-n=k$. In other words, the Beurling spectrum of an operator $A$ consists of the numbers of the non-zero diagonals of its matrix.
\eex

The module structures not conforming to the previous example are rarely considered in the literature. It is, however, an unnecessary and, at times, harmful restriction to consider only such kind of representations. In particular, within the framework of the above example one cannot consider more general Banach algebras than $B(\X)$ (see Remark \ref{ba}). The following structure, although very similar, also lies beyond the framework. 

\bex
For simplicity, we present this example in case $\GG = \RR^{2d}$ and leave an obvious extension
to more general LCA-groups to the reader. Consider the translation and modulation representations $T$ and $M$  defined by \eqref{trans1} and \eqref{mod1} on two
appropriate Banach spaces $\X$ and $\Y$ of functions on $\RR^d$. In this case, the representations
 satisfying our assumptions can be defined as 
\begeq\label{Weyl}
\T_{\X\Y}(t,s)A  =  M(t)T(s)AT(-s)M(-t), \ s,t\in\RR^d, \ A\in L(\X,\Y),\eq
and similarly for $\T_{\Y\X}$.
\eex

The following lemma is crucial for understanding the properties of the memory, i.e. Beurling spectra, of linear operators. Although a close analog of this result appears, for example, in \cite{BK05, BR75} we feel compelled to present a proof here. 

\bl\label{spprodf}
Let $A\in L(\X,\Y)$, $B\in L(\Y,\Z)$ and  $(L(\X,\Y), \T_{\X\Y})$, $(L(\Y,\Z), \T_{\Y\Z})$, and 
$(L(\X,\Z), \T_{\X\Z})$ be non-degenerate Banach modules over the algebras $L_{\nu_1}(\GG)$, $L_{\nu_2}(\GG)$, and $L_{\nu_3}(\GG)$, respectively, where $\nu_i$, $i=1,2,3$, are non-quazi-analytic weights. Assume also that \eqref{couple} is satisfied and that $A\in \overline{\Omega(A)}$ and
$B\in \overline{\Omega(B)}$ (see Definition \ref{modorbit}).
 Then
\[\L(BA,\T_{\X\Z})\subseteq\overline{\L(A,\T_{\X\Y})+\L(B,\T_{\Y\Z})}.\]
\el
\bpf
Assume that $\g\notin \overline{\L(A,\T_{\X\Y})+\L(B,\T_{\Y\Z})}=:\Delta$, and let
$\mathcal{U}$ and $\mathcal{V}$ be two neighborhoods of $\L(A,\T_{\X\Y})$
and $\L(B,\T_{\Y\Z})$, respectively, such that $\g\notin\mathcal U+\mathcal V$.
 Let also
$\phi\in L_{\nu_1}(\GG)$, $\psi\in L_{\nu_2}(\GG)$, and $f\in L_{\nu_3}(\GG)$ be such that $\supp \hat\phi\in \mathcal U$, $\supp\hat\psi\in\mathcal V$,  $\hat f(\g) \ne 0$, and $\supp \hat f \cap \overline{\supp \hat\phi +\supp\hat\psi}=\emptyset$. 
Such a function
$f$ exists because the weight $\nu_3$ is non-quazianalytic \cite{D56, LMF73}. 
In view of Lemma \ref{sprop}(vi) and \eqref{goodf} it is enough to prove that
$\L(B_1A_1,\T_{\X\Z})\subseteq\Delta$,
where $A_1 = \phi A$  and $B_1 = \psi B$.  
Then, omitting the indices of the representations, we get
\[
\bs
& \T(f)(B_1A_1)  = \int_\GG f(g)\T(-g)(B_1A_1)dg \\ &= \int_\GG f(g)(\T(-g)(\psi B))(\T(-g)(\phi A))dg \\
& = \int_\GG\int_\GG\int_\GG f(g)\phi(g_1)\psi(g_2)(\T(-g_2-g)B)(\T(-g_1-g)A)dg_1dg_2dg \\
& = \int_\GG\int_\GG\int_\GG f(g)\phi(g_1-g)\psi(g_2-g)(\T(-g_2)B)(\T(-g_1)A)dgdg_1dg_2 \\
& = \int_\GG\int_\GG F(g_1,g_2)(\T(-g_2)B)(\T(-g_1)A)dg_1dg_2.
\end{split}
\]
Observe that by our assumptions $\hat F(\g_1,\g_2) = \hat f(\g_1+\g_2)\hat\phi(\g_1)\hat\psi(\g_2) = 0$, and, hence, $\T(f)(B_1A_1) = 0$. Therefore, $\g\notin \Delta$ and the lemma is proved.
\epf

Lemma \ref{goodx} allows us to formulate the following  special case of the above result.

\bl\label{spprod}
Let $A\in L(\X,\Y)$, $B\in L(\Y,\Z)$ and  $(L(\X,\Y), \T_{\X\Y})$, $(L(\Y,\Z), \T_{\Y\Z})$, and 
$(L(\X,\Z), \T_{\X\Z})$ be non-degenerate Banach modules over the algebras $L_{\nu_1}(\GG)$, $L_{\nu_2}(\GG)$, and $L_{\nu_3}(\GG)$, respectively, where $\nu_i$, $i=1,2,3$, are non-quazi-analytic weights. Assume also that the above three algebras contain b.a.i., \eqref{couple} is satisfied, and that
$A$ and $B$ are $\T_{\X\Y}$- and $\T_{\Y\Z}$-continuous, respectively. Then
\[\L(BA,\T_{\X\Z})\subseteq\overline{\L(A,\T_{\X\Y})+\L(B,\T_{\Y\Z})}.\]
\el

The following corollary  appears, for example, in \cite{BK05, BR75}. 

\bc\label{spprodv}
Let $A \in (L(\X,\Y),\T_{\X\Y})_c$, where $\T_{\X\Y}$ is defined as in Example \ref{typop} and $x\in(\X,\T_\X)$. Then
\[\L(Ax,\T_{\Y})\subseteq\overline{\L(A,\T_{\X\Y})+\L(x,\T_{\X})}.\]
\ec

To help the reader, we note that the above is a generalization of the fact that a product of a one-diagonal matrix with the $m$-th non-zero diagonal and a one-diagonal matrix with the
$n$-th non-zero diagonal is either $0$ or another one-diagonal matrix with the $(m+n)$-th non-zero diagonal.

\brem\label{ba}
Analogous results hold if $A$, $B$ belong to a Banach algebra $\cB$, which has a non-degenerate $L_\nu(\GG)$-Banach module with respect to a strongly continuous representation $\T$ satisfying an obvious analog of \eqref{couple}, i.e. when $\{\T(g),\ g\in\GG\}$ is a group of automorphisms of the algebra $\cB$. In fact, most of the results about linear operators proved in this paper remain valid in the case of Banach algebras with a group of automorphisms and the proofs apply nearly verbatim. We shall discuss several examples in the last section of the paper.
\erem

We conclude this section with the following obvious, but useful, observation.

\bp\label{cinv}
Let $A \in (L(\X,\Y),\T_{\X\Y})_c$, where $\T_{\X\Y}$ is defined as in Example \ref{typop}.  Then
$A\X_c\subseteq\X_c$. If, moreover, $A$ is invertible then $A^{-1}\in(L(\Y,\X),\T_{\Y\X})_c$.
\ep 

\section{Memory decay for inverses of operators with compact Beurling spectrum}\label{Wcomp}

We begin with the simplest case of $\G = \RR$ and $\nu \equiv 1$. We consider an
invertible operator $A\in L(\X,\Y)$ with $\L(A) =\L(A,\T_{\X\Y})$ compact. We let $B = A^{-1}\in L(\Y,\X)$
and study the memory decay of $B$ with respect to the representation  $\T_{\Y\X}$. Because
$\nu \equiv 1$ the representations $\T_{\X\Y}$ and $\T_{\Y\X}$ are assumed to be isometric.

\bt\label{wcompact}
Assume $A \in L(\X,\Y)$ is invertible, $\L(A,\T_{\X\Y})\subset[-a,a]$,  and $B = A^{-1}$. Then there is $\alpha > 0$ such that  $B$ has exponential spectral decay of type $\alpha$, that is
$B\in D(\textbf{\emph{h}}_\a, \T_{\Y\X})$. Moreover, if $\phi_N \in L_1(\RR)$, $N\in\NN$, are
such that $\norm{\phi_N}=1$, 
$\hat\phi_N$ is continuous,  $\supp\hat\phi_N\subset [-N,N]$, and $\phi_{N,n}(x) = \phi_N(x)e^{ i Nnx}$, then
\begeq\label{infest}
\|\phi_{N,n} B\|\le \inf_{\a\in(0,\,  
a^{-1}\ln(1+\ae^{-1}(A))}\frac{e^{\alpha N(1-|n|))}
\|B\|}{1-(e^{\a a}-1)\aee(A)}, \ |n|>1,
\eq
where $\aee(A) = \|A\|\|B\|$ is the condition number of $A$.
\et

\bpf
Observe that since $\L(A)\subset [-a,a]$ Corollary \ref{Bineq}  implies that for every $z\in \CC$ we have $A\in D(e_z)$, where $e_z(\l) = e^{iz\l}$, and $A_\T: z\mapsto e_z\diamond A$ is an entire function.
Since invertibility is stable under small perturbations and the representation $\T$ is isometric, we get that $A_\T(z)$ is invertible for $z\in\CC_{-\d,\d}$ for some $\d > 0$. Clearly, then $(e_z\diamond A)^{-1}$ is a holomorphic extension of $B_\T$ and Lemma \ref{2side} implies $B\in D(\ee)$ for $0<\a<\d$. Obtaining the estimate \eqref{infest}, however, requires more work.

 For $z = t + i\a$,
$\a>0$,
let us estimate the norm $\|e_{z}\diamond A - \T_{\X\Y}(t)A\| = \|\textbf{e}_{-\a}\diamond A - A\|$. Consider
the $4a$-periodic function  $h(\l) = \eta(\l)-1$, where $\eta$ is defined by \eqref{curvehat}. We have shown in \eqref{fcoef} that
the Fourier coefficients $c_n$, $n\in\ZZ\backslash\{0\}$, of the periodic even function $\eta(\cdot - a)$ are positive. Clearly, they coincide with non-zero Fourier coefficients of $h(\cdot - a)$.  The $0$-th coefficient is also positive since
$c_0 = \frac4\a(\sinh \a a-\a a)>0$. Hence, by Proposition \ref{inmeas}, $\|\chT_{\X\Y}(h)\| \le 
e^{\a a}-1$.
Moreover, $h \equiv \textbf{e}_{-\a} -1\equiv e_{i\a}-1$
on  $[-a,a]\supseteq\L(A)$  and, hence, by Lemma \ref{sssprop}
\[\|\textbf{e}_{-\a}\diamond A - A\| = \|\chT_{\X\Y}(h) A\|
\le\|\chT_{\X\Y}(h)\|\|A\| = (e^{\a a}-1)\|A\|.\]
Using the same argument for the representation $\widetilde{\T}_{\X\Y}(t) = \T_{\X\Y}(-t)$,
we get 
\[\|\chT(\textbf{e}_{\a}) A - A\| = \|\textbf{e}_{\a}\diamond A - A\| \le (e^{\a a}-1)\|A\|,\ \a>0.\]

Therefore, if $|\a| <  a^{-1}\ln(1+\aee^{-1}(A))$, we have that
\[(\textbf{e}_{\a} \diamond A)^{-1} = B(I+(\textbf{e}_{\a} \diamond A-A)B)^{-1} = B\sum_{n=0}^\infty (-1)^n ((\textbf{e}_{\a} \diamond A-A)B)^n,\]
$B_\T(z): = (e_z \diamond A)^{-1} = \T_{\Y\X}(t)(\textbf{e}_{\mp\a} \diamond A)^{-1}$, $z=t\pm i\a$, 
$\a > 0$, satisfies
\[\|B_\T(z)\| \le \frac{\|B\|}{1-(e^{\a a}-1)\aee(A)}, \]
 and the estimate \eqref{infest} follows from \eqref{expest}.
\epf

Next, we extend the above result to the case when $\G = \RR^d$.

\bt
Assume $A \in L(\X,\Y)$ is invertible  and $B = A^{-1}$. 
Assume also that $a = (a_1,\ldots, a_d)\in\RR^d_+$ and $\L(A,\T_{\X\Y})\subset[-a_1,a_1]\times\ldots\times[-a_d,a_d]$. Then there is $\alpha = (\a_1,\ldots,\a_d) \in\RR^d_+$ such that the memory of $B$ has exponential decay of type $\alpha$, that is
\emph{$B\in D(\textbf{{h}}_\a, \T_{\Y\X})$}. Moreover, if $\phi^d_{N,n} \in L_1(\RR^d)$, $N\in\NN$, $n\in\ZZ^d$, are defined via \eqref{triangled}, then
\begeq\label{infestd}
\|\phi_{N,n} B\|\le \inf_{ \a  
}\frac{e^{ NM_n(\a)}
\|B\|}{1-\left(e^{\a\cdot a}-1\right)\aee(A)}, \ n\in \ZZ^d,
\eq
where 
$M_n(\a) = \sum\limits_{k:\ n_k\ne 0} \a_k(1-|n_k|)$ and the infimum is taken over all $\a\in\RR^d$ such that the denominator is positive.
\et

\bpf
The proof of the one-dimensional case can be applied almost verbatim if the periodic function $\eta$ is replaced with
$\eta^d(\l) = \prod_{k=1}^d \eta(\l_k)$, $\l = (\l_1,\ldots, \l_d)\in\RR^d$.
Remark \ref{tenp} and Lemmas \ref{ssprop} and  \ref{sssprop}, however, cannot be used in this case. To overcome this minor obstacle one uses Lemma \ref{sprop} to prove the estimates with $a+\varepsilon$, $\varepsilon\in\RR^d_+$, in place of $a$. The desired estimates then follow since $\varepsilon\in\RR^d_+$ is arbitrary.
We omit the remaining details for the sake of brevity.
\epf

\section{Memory decay of inverses to causal operators and operators with exponential memory decay}\label{Wcaus}
In this section we extend the above results to operators with exponential spectral decay as well as to causal operators and operators with one-sided exponential spectral decay. As before the proof is based on the holomorphic extensions of the function $t\mapsto \T(t)A$. We begin, however, with establishing the algebraic properties of the classes of operators considered in this section.

\bl\label{causalg}
Assume $A\in\C(L(\Y,\Z))$ and $B\in\C(L(\X,\Y))$. Then $AB\in\C(L(\X,\Z))$.
\el
\bpf
Follows from Lemma \ref{spprod}.
\epf

\bl\label{prodexp1}
Let $\a\in\RR^d$. Assume $A\in L(\Y,\Z)_c\cap D(\emph{\textbf{e}}_{\a}, \T_{\Y\Z})$ and 
$B\in L(\X,\Y)_c\cap D(\emph{\textbf{e}}_{\a}, \T_{\X\Y})$. Then $AB\in L(\X,\Z)_c\cap D(
\emph{\textbf{e}}_{\a}, \T_{\X\Z})$ and
\begeq\label{hom1side}
\emph{\textbf{e}}_{\a}\diamond(AB) = (\emph{\textbf{e}}_{\a}\diamond A)( \emph{\textbf{e}}_{\a}\diamond B).  
\eq
\el
\bpf
The proof follows by uniqueness of the holomorphic extension provided by Proposition \ref{1expestprop}.
\epf

\bt\label{1sidexpwin}
Let $\a\in\RR^d_\pm$ and $A\in L(\X,\Y)_c\cap D(\emph{\textbf{e}}_{\a})$. Assume
$B = A^{-1} \in L(\Y,\X)$. Then there is $\beta\in \RR^d_\pm$ such that 
$B\in L(\Y,\X)_c\cap D(\emph{\textbf{e}}_{\beta})$.
\et

\bpf
All we need to do is repeat the first part of the proof of Theorem \ref{wcompact}. Indeed, if $A_\T$ is
the holomorphic extension provided by Proposition \ref{1expestprop}, it is continuous in the uniform operator topology in
$\overline{\CC_{\mp\d}}$ for some $\d\in\RR^d_+$. Hence, because  the representation $\T$ is bounded and invertibility is
stable under small perturbations, there exists $\beta \in \RR^d_+$ such that $A_\T(z)$ is invertible for  $z\in\overline{\CC_{\mp\beta}}$. It remains only to apply Proposition \ref{1expestprop} once again.
\epf

The following two results are immediate consequences of Theorem \ref{1sidexpwin}. A version of the first of them was originally announced in \cite{BK06}.

\bt\label{causwin}
Assume $A\in L(\X,\Y)_c\cap \C_\pm(L(\X,\Y))$ and $B = A^{-1}\in L(\Y,\X)$. Then there is $\a\in \RR^d_+$ such that 
$B\in L(\Y,\X)_c\cap D(\emph{\textbf{e}}_{\mp\a})$.
\et

\bt\label{invexp}
Let $\a\in\RR^d_+$ and $A\in L(\X,\Y)_c\cap D(\emph{\textbf{h}}_{\a})$. Assume
$B = A^{-1} \in L(\Y,\X)$. Then there is $\beta\in \RR^d_+$ such that 
$B\in L(\Y,\X)_c\cap D(\emph{\textbf{h}}_{\beta})$.
\et

\brem
Guided by the fact that the inverse of a triangular matrix is also triangular, one might think that the inverse to a causal operator should always be causal. The bilateral shift operator shows, however, that the inverse to a causal operator may, in fact, be anticausal. Theorem \ref{causwin}, on the other hand, ensures that the anticausal part of the inverse always has exponential memory decay. We refer to \cite{BK05} for different criteria of causal invertibility.
\erem

\brem
Given an invertible operator $A\in D(\eet)$ we cannot obtain estimates for the memory decay
of $A^{-1}$ without assuming more than $A\in L(\X,\Y)_c$. It is possible, however, to use the approach of Theorem \ref{wcompact} to  obtain such estimates if we can quantify the convergence of
$\T(t)A\to A$ as $t\to 0$. In particular, we can get the estimates if $A$ is in one of the H\"older-Zygmund classes studied in \cite{GK10}.
\erem

\section{Extension of the classical Wiener's Lemma}\label{Wclass}

We begin by establishing the algebraic property of the Wiener Class $\W = \W(L(\X,\Y))$ of operators in $L(\X,\Y)$. Recall from Section \ref{wiencl} that the norms of operators $A\in\W$ are given by 
\begeq\label{wiensum}
\|A\|_{1,N} = 5^d \sum_{n\in\ZZ^d} \|\phi^d_{N,n} A\|= 5^d \sum_{n\in\ZZ^d} \|\T_{\X\Y}(\phi^d_{N,n}) A\|,
\ N\in\NN,
\eq
where we used the functions $\phi_N^d$ defined by
\eqref{triangled} and their translates $\phi^d_{N,n}$.

\bl\label{algprop}
If 
$B\in\W(L(\X,\Y))$ and $A\in\W(L(\Y,\Z))$ then $AB\in\W(L(\X,\Z))$ and
\[\|AB\|_{1,N}\le\|A\|_{1,N}\|B\|_{1,N},\]
in particular, $\W(B(\X))$ is a Banach algebra.
\el
\bpf
We deduce the desired property from
\begin{equation*}
\begin{split}
\|AB\|_{1,N} & = 5^d \sum_{n\in\ZZ^d} \|\phi^d_{N,n} (AB)\| \\ &= 
5^d \sum_{n\in\ZZ^d} \left\|\phi^d_{N,n} \left(\left(\sum_{m\in\ZZ^d} \phi^d_{N,m}A\right)
\left(\sum_{k\in\ZZ^d} \phi^d_{N,k}B\right)\right)\right\| \\ &
\le 5^{2d} \sum_{m\in\ZZ^d} \|\phi^d_{N,m} A\|\sum_{k\in\ZZ^d} \|\phi^d_{N,k} B\| = \|A\|_{1,N}\|B\|_{1,N},
\end{split}
\end{equation*}
where the inequality is true because for every choice of $m,n\in\ZZ^d$ there are at most $5^d$ different numbers $k\in\ZZ^d$ such that the set $\supp \hat\phi^d_{N,k}\cap \L((\phi^d_{N,m}A)( \phi^d_{N,n}B),\T_{\X\Z})$ has non-empty interior. We used Lemmae \ref{ssprop}, \ref{spprod} and Remark \ref{tenp}.
\epf

\bt\label{maint}
Let $A\in \W(L(\X,\Y))$ be invertible with $B = A^{-1} \in L(\Y,\X)$. Then $B\in \W(L(\Y,\X))$ and
\begeq\label{wienest}
\|B\|_{1,1} \le (2N+1)^d\|B\|_{1, N} \le\frac{\|B\|}{\eps_A}\left(2\psi_A\left(\frac{\eps_A}{\|B\|}\right)+1\right)^d,
\eq
where $N = \psi_A(\frac{\eps_A}{\|B\|})$, $\eps_A = 5^{-d}\left(16\left(\frac{2\d_A - 1}{\d_A-1}\right)^d-12\right)^{-1}$, 
$\d_A = \left(\frac{4\ae(A)+3}{4\ae(A)+2}\right)^{\frac1{3d}}$, $\aee(A) = \|A\|\|A^{-1}\|$, and $\psi_A: \RR_+\to \NN$ is defined by
\[
\psi_A\left(t\right) = \min \left\{K\in\NN:\ \norm{A - \sum_{|k|_\infty \le 2} \phi^d_{K,k}A}_{1,K} \le t\right\},
\]
where $\phi^d_{K,k}$ is given by \eqref{triangled}.
\et
\bpf
We begin by fixing some $N\in\NN$, which is to be determined later, and using the functions 
$\phi_{N,n}^d$ defined via \eqref{triangled} to represent $A = C+D$, $C =  \sum\limits_{|n|_\infty \le 2} \phi^d_{N,n}A$,
in such a way that 
$\|B\|\|D\|_{1,N} \le\frac12$ and, hence, also $\|B\|\|D\| \le\frac12$. The first inequality is true for some $N\in\NN$ because
$A\in \W(L(\X,\Y))$. In this case, $\|C\|\le \|A\|+\frac1{2\|B\|}$, $C$ is invertible, and 
$L = C^{-1}=(A-D)^{-1} = B(I-DB)^{-1} = B\sum_{n=0}^\infty (DB)^n$ satisfies
\begeq
\|L\|\le \frac{\|B\|}{1-\|B\|\|D\|} \le 2 \|B\|.
\eq
Using \eqref{infestd} for the operators $C$ and $L = C^{-1}$ with $a= (3N,\ldots,3N)\in\RR^d$, 
we get 
\begeq
\bs
\|\phi^d_{N,n}L\| &\le \inf_{\a} \frac{e^{\a NM_n}\|L\|}{1-(e^{3\a d N}-1)\aee(C)} \\ &\le
 \inf_{\a} \frac{e^{\a NM_n}\|B\|}{1-\|B\|\|D\|-(e^{3\a dN}-1)(\aee(A)+\frac12)},
 \\ &\le \inf_{\a} \frac{e^{\a NM_n}\|B\|}{\frac12-(e^{3\a dN}-1)(\aee(A)+\frac12)},
\end{split}
\eq
where $M_n = \sum\limits_{k=1}^d 1-|n_k|$, $n=(n_1,\ldots, n_d)\in\ZZ^d$, $\ n_k\ne 0$, and the infima are taken over all $\a\in\RR_+$ such that the corresponding denominator is positive.

Next, we choose $\a$ such that
\[\frac12-(e^{3\a dN}-1)(\aee(A)+\frac12) = \frac14.\]
Then $e^{\a N} = \left(\frac{4\ae(A)+3}{4\ae(A)+2}\right)^{\frac1{3d}}=:\d_A > 1$ and, therefore,
\begeq
\|\phi^d_{N,n}L\| \le 4\d_A^{M_n}\|B\|, \ n\neq 0, \mbox{ and } \|\phi^d_{N,0} L\|\le\|L\|\le 2\|B\|.
\eq
Summing over $n\in\ZZ^d$ and multiplying by $5^d$, we get
\begeq
\|L\|_{1,N} \le 5^d\left(8\left(\frac{2\d_A - 1}{\d_A-1}\right)^d-6\right)\|B\| = \frac1{2\eps_A}\|B\|.
\eq
Finally, we specify some $N\in\NN$ such that
\[
\|L\|_{1,N}\|D\|_{1,N} \le \frac1{2\eps_A}\|B\|\|D\|_{1,N} \le\frac12,
\]
that is, we let 
$N = \psi_A\left( \frac{\eps_A}{\|B\|}\right)$. Then the series in $B = (C + D)^{-1} = L(I+DL)^{-1} = L\sum_{n=0}^\infty (-1)^{n}(DL)^n$ converges in $\W(B(\Y))$, and,
using
\[\|B\|_{1,N} \le\frac{\|L\|_{1,N}}{1-\|L\|_{1,N}\|D\|_{1,N} }\le \frac{\|B\|}{\eps_A},\]
we obtain the estimate \eqref{wienest}. 
\epf

\bp\label{invarw}
Assume $A\in\W(L(\X,\Y),\T_{\X\Y})$, where the representation $\T_{\X\Y}: \RR^d\to B(L(\X,\Y))$ is defined as in Example \ref{typop} by $\T_{\X\Y}(t)A = \T_\Y(t)A\T_\X(-t)$. Then 
\[A\W(\X,\T_\X) \subseteq \W(\Y,\T_\Y).\]
Moreover, if $A$ is invertible then it is also an isomorphism between   $\W(\X,\T_\X)$ and $\W(\Y,\T_\Y)$.
\ep

\bpf
Let $x\in \W(\X,\T_\X)$ and $y = Ax$. As in Lemma \ref{algprop}, we get
\begin{equation*}
\bs
\sum_{k\in\ZZ^d} \|\phi^d_{N,k}y \| &= 
\sum_{k\in\ZZ^d} \left\|\phi^d_{N,k} \left(\left(\sum_{m\in\ZZ^d} \phi^d_{N,m}A\right)
\left(\sum_{n\in\ZZ^d} \phi^d_{N,n}x\right)\right)\right\|  \\
& \le 5^{d} \sum_{m\in\ZZ^d} \|\phi^d_{N,m} A\|\sum_{n\in\ZZ^d} \|\phi^d_{N,n} x\| = \|A\|_{1,N}\|x\|_{1,N},
\end{split}
\end{equation*}
where we used Corollary \ref{spprodv}.
\epf

\brem
Various corollaries of the above proposition play a crucial role in the study of canonical duals of localized frames (see, e.g. \cite{BCHL06I, F09, G04, KO08} and references therein).
\erem

\section{Inverse Closedness of the Beurling Class}\label{Wbeur}
The Beurling Class $\B = \B(\X,\Y)$ of operators in $L(\X,\Y)$ has the norm
\begeq\label{bnormop}
\|A\|_{\B}=  \sum_{k\in\ZZ^d} \max_{|n|_\infty\ge |k|_\infty} \|\phi^d_{1,n}A\| = 
 \sum_{k\in\ZZ^d} \max_{|n|_\infty\ge |k|_\infty}\|\T_{\X\Y}(\phi^d_{1,n})A\|,
\eq
where $|k|_\infty = \max\{|k_1|, |k_2|,\ldots, |k_d|\}$, $k = (k_1, k_2,\ldots, k_d)\in\RR^d$. 

We prove the inverse closedness of this class by adapting a similar argument in \cite{S10ca} to our general setting and applying the technique of the previous section.

We begin with  the following
very strong algebraic property.

\bl
If 
$B\in\B(\X,\Y)$ and $A\in\B(\Y,\Z)$ then $AB\in\B(\X,\Z)$ and
\begeq\label{brand}
\|AB\|_{\B}\le 5^d\|A\|_{1,1}\|B\|_{\B}+2^d\|A\|_{\B}\|B\|_{1,1},
\eq
in particular, $\B(\X)=\B(\X,\X)$ is a Banach algebra with respect to an equivalent norm.
\el
\bpf
Using the absolute convergence of the series below, we have 
\begin{equation*}
\begin{split}
\|\phi_{1,n}(AB)\| & =\left\|\phi_{1,n}\left(\left(\sum_{i\in\ZZ^d} \phi_{1,i}A\right)\left(\sum_{j\in\ZZ^d}
\phi_{1,j}B\right)\right)\right\| \\
& = \left\|\sum_{i\in\ZZ^d}\sum_{j\in\ZZ^d}\phi_{1,n}\left(\left(\phi_{1,i}A\right)\left(
\phi_{1,j-i}B\right)\right)\right\| \\
&
\le  \left\|\sum_{|i|_\infty \ge\frac12|n|_\infty}\sum_{j\in\ZZ^d}\phi_{1,n}\left(\left(\phi_{1,i}A\right)\left(
\phi_{1,j-i}B\right)\right)\right\| \\ &+ \left\|\sum_{|i|_\infty <\frac12|n|_\infty}\sum_{j\in\ZZ^d}\phi_{1,n}\left(\left(\phi_{1,i}A\right)\left(\phi_{1,j-i}B\right)\right)\right\| 
\end{split}
\end{equation*}
From Lemma \ref{spprod} we deduce that when $|j-n|_\infty >2$ the terms in both series vanish. Hence,

\begin{equation*}
\bs
\|\phi_{1,n}(AB)\|  &\le \sum_{|i|_\infty \ge\frac12|n|_\infty}\sum_{|j-n|_\infty \le 2}\|\phi_{1,i}A\| \| \phi_{1,j-i}B\| \\
& +\sum_{|j-n|_\infty \le 2}\sum_{|j-i|_\infty <\frac12|n|_\infty}\|\phi_{1,j-i}A\| \|\phi_{1,i}B\| 
\\
& \le \sum_{|i|_\infty \ge\frac12|n|_\infty}\sum_{|j-n|_\infty \le 2}\|\phi_{1,i}A\| \| \phi_{1,j-i}B\| \\
& +\sum_{|i|_\infty \ge\frac12|n|_\infty-1}\sum_{|j-n|_\infty \le 2}\|\phi_{1,j-i}A\| \|\phi_{1,i}B\|.
\end{split}
\end{equation*}
Therefore,
\begin{equation*}
\bs
\sum_{k\in\ZZ^d} \max_{|n|_\infty\ge |k|_\infty} & \|\phi_{1,n}(AB)\| \\ &\le 
\sum_{k\in\ZZ^d} \max_{|n|_\infty\ge |k|_\infty} \sum_{|i|_\infty \ge\frac12|n|_\infty}\sum_{|j-n|_\infty \le 2}\|\phi_{1,i}A\| \| \phi_{1,j-i}B\| \\
& +\sum_{k\in\ZZ^d} \max_{|n|_\infty\ge |k|_\infty}\sum_{|i|_\infty \ge\frac12|n|_\infty-1}\sum_{|j-n|_\infty \le 2}\|\phi_{1,j-i}A\| \|\phi_{1,i}B\| \\
&\le \sum_{k\in\ZZ^d}  \max_{|n|_\infty\ge \frac { |k|_\infty}2}\|\phi_{1,n}A\|\sum_{i\in\ZZ^d}5^d\|\phi_{1,i}B\|  \\
&+\sum_{k\in\ZZ^d}  \max_{|n|_\infty\ge \frac { |k|_\infty}2-1}\|\phi_{1,n}B\|\sum_{i\in\ZZ^d}5^d\|\phi_{1,i}A\|,
\end{split}
\end{equation*}
and the desired inequality follows. Since $\|A\|_\W\le 5^d\|A\|_\B$, the inequality implies that $\B(\X)$ is indeed a Banach algebra.
\epf

Next we use the standard Brandenburg trick \cite{B75} to show that the spectral radii in the Wiener and
Beurling algebras coincide, i.e., for all $A\in \B(\X)$
\begeq\label{radii}
\rho_\B(A) = \lim_{n\to\infty}\|A^n\|^{1/n}_{\B} = \lim_{n\to\infty}\|A^n\|^{1/n}_{\W} = \rho_\W(A).
\eq
Indeed, since $\|A\|_\W\le 5^d\|A\|_\B$, we have $\rho_\W(A)\le\rho_\B(A)$. On the other hand,
from \eqref{brand} we get $\|A^{2n}\|_\B \le (5^d+2^d)\|A^n\|_\B\|A^n\|_{1,1}$. Taking the $n$-th root and
passing to the limit we get $\rho_\B(A^2) = (\rho_\B(A))^2 \le \rho_\B(A)\rho_\W(A)$ which implies
\eqref{radii}.

At this point, Hulanicki's lemma \cite{H72} is typically used to finish the argument about inverse closedness \cite{G10, GK10, S10}. In our setting, however, this result is not applicable since we may be entirely outside of the realm of $*$-algebras. We overcome this obstacle by looking more closely into the proof of the result in the previous section.

\bt\label{absB}
Assume that $A\in\B(\X,\Y)$ and $B\in L(\Y,\X)$ is its inverse operator. Then $B\in\B(\Y,\X)$.
\et
\bpf
As in the proof of Theorem \ref{maint}, let $A = C + D$ where $\L(C)$ is compact and
\[B = C^{-1}(I+DC^{-1})^{-1} = C^{-1}\sum_{k=0}^\infty (-DC^{-1})^n,\]
where the Neumann series converges in $\W(\Y)$ because $\|DC^{-1}\|_\W < 1$. 
Theorem \ref{wcompact} implies $C^{-1}\in D(\ee)\subseteq \B$  and we can use \eqref{radii} to conclude that $\rho_\B(DC^{-1}) < 1$ and, hence, the Neumann series converges in $\B(\Y)$. 
\epf

The proof of the following result is analogous to that of proposition \ref{invarw} and is left to the reader.

\bp\label{invarb}
Assume $A\in\B(L(\X,\Y),\T_{\X\Y})$, where the representation $\T_{\X\Y}: \RR^d\to B(L(\X,\Y))$ is defined as in Example \ref{typop} by $\T_{\X\Y}(t)A = \T_\Y(t)A\T_\X(-t)$. Then 
\[A\B(\X,\T_\X) \subseteq \B(\Y,\T_\Y).\]
Moreover, if $A$ is invertible then it is also an isomorphism between   $\B(\X,\T_\X)$ and $\B(\Y,\T_\Y)$.
\ep

\section{Sobolev-type Algebras}\label{sobol}

Here we extend the previous results to the Sobolev-type classes $\W^{(m)}$ and $\B^{(m)}$.
Unlike the previous two sections, rather than deal with operators in $L(\X,\Y)$ we choose to work with a unital Banach algebra $\cB$ such that the representation $\T$ defines a group of isometric automorphisms of $\cB$. It should be clear how to modify the statements to get the results for $L(\X,\Y)$. The definitions of $\W^{(m)}$ and $\B^{(m)}$ were
given in Subsection \ref{sclass}; in what follows we are using the same notation. 

 In pursuing our extension, we follow the ideas in \cite{BR75, GK10}. We note, however, that since
our Wiener and Beurling classes are different from those considered in \cite{GK10}, the Sobolev-type classes are also different, and the results in this section do not follow immediately from \cite{GK10}.

Observe that since each $\T^{(k)}(t)$ defined in Remark \ref{tenp} is an automorphism of $\cB$, the definition of the infinitesimal generator $A_k$ implies
\begeq\label{deriv}
A_k(xy) = (A_k x)y + x (A_k y), \ x,y\in D(A_k),
\eq
i.e. each $A_k$, $k = 1,\ldots, d$, is a derivation \cite{BR75} on $\cB$. 
These equalities ensure that both $\W^{(m)}$ and $\B^{(m)}$ are Banach algebras with the norm
\[\|x\|_{\F^{(m)}} = \sum_{|\a|\le m} \|A^\a x\|_\F,\ x\in\F,\]
where $\F$ is either $\W$ or $\B$. We refer to these algebras as Sobolev-Wiener or Sobolev-Beurling algebras, respectively.

\bt\label{sobolt}
Let $\F^{(m)}$, $m\in\NN$, be either a Sobolev-Wiener or a Sobolev-Beurling algebra. Assume that
$x \in\F^{(m)}$ is invertible, i.e. $x^{-1}\in \cB$. Then $x^{-1}\in \F^{(m)}$, i.e. $\F^{(m)}$ is an inverse-closed subalgebra of $\cB$.
\et
\bpf
Observe that $A_k e = 0$ for the unit $e\in\cB$ and $D(A_k)$ is inverse closed since the function $t\mapsto\T^{(k)}(t)x^{-1} = (\T^{(k)}(t)x)^{-1}$ is differentiable when $x\in D(A_k)$ is invertible in $\cB$
(see also \cite{BR75, GK10}). Therefore, using equalities \eqref{deriv} we get
$0 = A_k(xx^{-1}) = (A_k x)x^{-1} + x(A_k x^{-1})$ and, hence,
\[A_k x^{-1} = -x^{-1}(A_k x) x^{-1}.\]
The result is now an immediate consequence of Theorems \ref{maint}, and \ref{absB}.
\epf

\bp\label{invarsvb}
Assume $C\in\F^{(m)}(L(\X,\Y),\T_{\X\Y})$, where the representation $\T_{\X\Y}: \RR^d\to B(L(\X,\Y))$ is defined as in Example \ref{typop} by $\T_{\X\Y}(t)C = \T_\Y(t)a\T_\X(-t)$. Then 
\[C\F^{(m)}(\X,\T_\X) \subseteq \F^{(m)}(\Y,\T_\Y).\]
Moreover, if $C$ is invertible then it is also an isomorphism between   $\F^{(m)}(\X,\T_\X)$ and $\F^{(m)}(\Y,\T_\Y)$.
\ep

\bpf
Let $A^{\X}$, $A^{\Y}$ and $A^{\X\Y}$ be the multigenerators of the representations  $\T_\X$, $T_\Y$ and $\T_{\X\Y}$, respectively. Then, using the relation between the representations, we obtain
\begeq
A_k^{\Y}(Cx) = (A_k^{\X\Y} C)x + C (A_k^\X x), \ x\in D(A^\X_k),
\eq
similarly to \eqref{deriv}. It remains to apply Propositions \ref{invarw} and \ref{invarb}.
\epf

\section{Examples}\label{examp}

\subsection{Matrices and resolutions of the identity}\label{matrex} 
We begin by revisiting the example we presented in the introduction and proceed to illustrate how general, in fact, it is or can be made.

Let us consider  operators in $L(\X,\Y)$ where $\X$ and $\Y$ are Banach modules with the structures associated with the representations $\T_\X$ and $\T_\Y$, respectively. Assume that the representations $\T_{\X\Y}$ and $\T_{\Y\X}$ are defined as in the example \ref{typop}.
Let $\s_k$, $k\in\ZZ^d$, be a partition of $\RR^d$ such that
\begin{itemize}
\item The spectral submodules $\X(\s_k,\T_\X)$, $\Y(\s_k,\T_\Y)$, $k\in\ZZ^d$, are complementable Banach subspaces, i.e.~for each $k\in\ZZ^d$ there exist two projections  $P_k\in B(\X)$ and $Q_k\in B(\Y)$ onto $\X(\s_k, \T)$ and $\Y(\s_k, \T)$ respectively;
\item Projections $P_k$, $Q_k$, $k\in\ZZ^d$, form  disjunctive resolutions of the identity $\cP$ and $\cQ$ in $\X$ and $\Y$ respectively;
\item $\cP$ and $\cQ$ satisfy \eqref{crp}.
\end{itemize}
Observe that the above conditions are typically satisfied in Hilbert spaces as well as in different kinds
of $L^p$ spaces. 

Let us now define the representations $U_\X$, $U_\Y$ by \eqref{exres} and $U_{\X\Y}$ by \eqref{typop1}. Now we can define various classes $\F(L(\X,\Y), U_{\X\Y})$ in terms of the matrix
$(Q_kAP_j)_{k,j\in\ZZ^d}$ of an operator $A\in L(\X,\Y)$ as it was done in \cite{B97Izv}. We leave it to 
the reader to obtain the specific formulas. We only notice that if 
$\s_k = k+[0,1)^d$, $k\in\ZZ^d$, computations similar  to those in the proof of Proposition \ref{preq}
show that $\F(L(\X,\Y),U_{\X\Y}) = \F(L(\X,\Y),\T_{\X\Y})$.

 Alternatively, let $\cB$ be a (complex) unital Banach algebra and $\T: \RR^d\to B(\cB)$ be a group
of bounded (isometric) automorphisms of $\cB$. Let $\s_k$, $k\in\ZZ^d$, be a partition of $\RR^d$ such that
\begin{itemize}
\item The spectral submodules $\cB(\s_k,\T)$, $k\in\ZZ^d$, are complementable Banach subspaces, i.e. for each $k\in\ZZ^d$ there exists a projection $P_k\in B(\cB)$ onto $\cB(\s_k, \T)$;
\item Projections $P_k$, $k\in\ZZ^d$, form a disjunctive resolution of the identity $\cP$;
\item $\cP$ satisfies \eqref{crp}.
\end{itemize}
Given the above conditions we can define a bounded $2\pi$-periodic representation $U_\cB: \RR^d\to B(\cB)$ similar to \eqref{exres}:
\[U_\cB(t)x = \sum_{k\in \ZZ^d} e^{int}P_kx, \ t\in\RR^d, x\in\cB.\]
Various classes $\F(\cB, U_\cB)$ are then easily defined in terms of $\|P_kx\|$, $k\in \ZZ^d$, as e.g.~ in Example \ref{ex1}. Again,
if $\s_k = k+[0,1)^d$, $k\in\ZZ^d$, computations similar  to those in the proof of Proposition \ref{preq}
show that $\F(\cB,U_\cB) = \F(\cB,\T)$.

Clearly, it is unreasonable to expect the conditions on the spectral submodules outlined in this example to hold outside of the realm of Hilbert and $L^p$ spaces. The following examples
are designed to show that the techniques of this paper still produce meaningful results when previously used methods fail.

\subsection{Continuous operator-valued functions}\label{1examp}
 Let $\cB$ be a (complex) unital Banach algebra. In this example we consider the Banach algebra $C_b = C_b(\RR^d, \cB)$ of all bounded
continuous $\cB$-valued functions with the norm $\|x\| = \sup_{s\in\RR^d} \|x(s)\|$, $x\in C_b$.
The algebra $C_b$ is equipped with a non-degenerate $L_1(\RR^d)$-module structure associated with the translation representation $\T = T: \RR^d\to B(C_b)$ similar to \eqref{trans1}:
\[(T(t)x)(s) = x(t+s),\ t,s\in\RR^d,\ x\in C_b.\]
In particular, for $f\in L_1(\RR^d)$ the module structure is given by the convolution
\[(T(f)x)(s) = (f*x)(s) = \int_{\RR^d} f(t)x(s-t)dt, \ x\in C_b.\]
The Beurling spectrum $\L(x,T)$ in this case coincides with $\supp \hat x$, where the Fourier transform is understood in the sense of tempered distributions. The class of $T$-continuous functions in this case is
the algebra $C_{ub}=C_{ub}(\RR^d,\cB)$ of all bounded uniformly continuous $\cB$-valued functions.

Consider the functions $\phi_{1,a}^d$ defined by \eqref{triangled}. The following subclasses of $C_{ub}$ (see Remark \ref{wc}) correspond to the classes considered in the previous sections. The Wiener class is
\[\W = \W(C_b)= \W(C_b, T) = \{x\in C_b:\ \int_{\RR^d} \|\phi_{1,a}^d*x\| da <\infty\},\]
the class of functions with exponential spectral decay is
\[
\bs
\W_{exp}  
= \bigcup_{\a\in\RR^d_+}D(\ee) =  \ & \{x\in C_b:  \|\phi_{1,n}^d*x\| \le M e^{-\varepsilon|n|}, n\in\ZZ,\\ &
\mbox{ for some positive } M=M(x) \mbox{ and } \varepsilon = \varepsilon(x)\},
\end{split}
\]
the classes of functions with one-sided exponential spectral decay are
\[
\bs
\W_{exp}^\pm  
=  \bigcup_{\a\in\RR^d_\pm}D(\eet) =  \ & \{x\in C_b:  \|\phi_{1,n}^d*x\| \le M e^{-\varepsilon|n|}, \pm n \in \NN,\\ &
\mbox{ for some positive } M=M(x) \mbox{ and } \varepsilon = \varepsilon(x)\},
\end{split}
\]
and the Beurling class is
\[\B = \B(C_b) = \{x\in C_b:\ \sum_{k\in\ZZ^d} \max_{|n|_\infty\ge |k|_\infty} \|\phi^d_{1,n} * x\| < \infty\}.\]
A straightforward computation shows that the Sobolev-type algebras in this case are
\[
\bs
\W^{(m)} & = \W^{(m)}(C_b) = \{x\in C^m_b:\ \sum_{|\a|\le m}\int_{\RR^d} \|\phi_{1,a}^d*x^{(\a)}\| da <\infty\} \\
& = \{x\in C_b:\ \int_{\RR^d} \|\phi_{1,a}^d*x\|(1+|a|)^m da <\infty\},
 \end{split}
 \]
 and
\[\B^{(m)} = \B^{(m)}(C_b) = \{x\in C^m_b:\ \sum_{|\a|\le m}\sum_{k\in\ZZ^d}^\infty \max_{|n|_\infty\ge |k|_\infty} \|\phi^d_{1,n} * x^{(\a)}\| < \infty\},\]
where $C_b^m$ is the space of functions with continuous bounded derivatives up to order $m\in\NN$ and 
$x^{(\a)}$ is the (higher order) partial derivative corresponding to the index $\a\in\RR^d$. 

In the following theorem we denote by $\F(C_b)$ one of the above classes.

\bt\label{cb}
Assume that $x \in \F(C_b)$ is such that $x(s)$ is invertible in $\cB$ for all $s\in \RR^d$ and the function $y(\cdot) = (x(\cdot))^{-1}$ is bounded. Then $y\in\F(C_b)$.
\et

\bpf
The result follows immediately from Theorems \ref{1sidexpwin}, \ref{invexp}, \ref{maint}, \ref{absB}, \ref{sobolt} and Remark \ref{ba}.
\epf


\brem
In case $\cB = \CC$ and $d=1$, the classes $\W_1$ and $\W_{exp}$ were considered in \cite{B97Sib}. The following two families of  classes indexed by $q>1$ were also shown to be inverse closed there:
\[\W_q = \{x\in C_b(\RR):\ \int_{\RR^d} \|\phi_{1,a}^d*x\|(1+|a|)^q da <\infty\} \mbox{ and }\]
\[\W_{q,0} = \{x\in C_b(\RR):\  \lim_{|a|\to\infty}\|\phi_{1,a}^d*x\| |a|^q =0\}.\]
\erem

\brem
Observe that if $x\in C_b(\RR)$ is a periodic function, $x(t) \simeq \sum_{n\in\ZZ} c_ne^{int}$, then $x\in \W$ if and only if it has summable Fourier coefficients. Indeed, in this case, we have
$c_ne^{int} = (\phi_{1,n}*x)(t)$, $n\in\ZZ$. Hence, the classical Wiener's lemma \cite{W32} follows from Theorem \ref{cb}. 
\erem

\brem
In general, if  $x\in C_b(\RR^d)$, then $|(\phi_{1,a}^d*x)(t)|$ coincides with the absolute value of the Short-time Fourier 
transform (STFT) $|(V_{\phi}x)(t,a)|$ of the function $x$ with respect to the window $\phi$ \cite{G01}. Hence, $\W(C_b)$ coincides with the modulation space $M^{\infty,1}(\RR^d)$ and the above theorem implies that if a function $x\in M^{\infty,1}(\RR^d)$ is bounded away from $0$ then
$1/x \in M^{\infty,1}(\RR^d)$ \cite{S95}. In view of the described relation with the STFT, one may call 
the map $x \mapsto (\mathcal{V}_\phi x)(a) =\T(\phi_{1,a}x)$ the \emph{abstract STFT} of the element $x$ with respect to a window $\phi$ and a representation $\T$. It is also natural to ask if the methods in this paper can be used to recover the results in \cite{S94, S95, G06}, that is to show that
$M^{\infty,1}(\RR^d)$, also known as the Sj\"ostrand's class, is an inverse closed subalgebra when the product is defined via twisted convolution as opposed to the ordinary convolution. We believe, that the answer is positive but non-trivial and, hence, deserves a separate paper.
\erem

\brem
We have presented Theorem \ref{cb} as a special case of Theorems
\ref{1sidexpwin}, \ref{invexp}, \ref{maint}, \ref{absB}, and \ref{sobolt}.
 In fact, this case is generic when the algebras are considered.
Indeed, let $\cB$ be a unital Banach  algebra with a strongly continuous group of automorphisms 
$\T:\RR^d\to B(\cB)$. 
We can regard $\cB$ as an inverse closed subalgebra of $C_b(\RR^d,\cB)$ by identifying $x\in\cB$ with $x_\T \in C_b$,
$x_\T(t) = \T(t)x$. Moreover, we then have $T(s)x_\T(t) = x_\T(t+s) = \T(t+s)x$ and, hence,
\[(T(f)x_\T)(t) = \int_{\RR^d} f(s)\T(t-s)xds = \T(t)(fx) = (fx)_\T(t),\ f\in L_1(\RR^d).\]
Therefore, $\L(x,\T) = \L(x_\T,T)$ and $x\in \F(\cB)$ if and only if $x_\T\in\F(C_b)$, where $\F$ denotes one of the classes in Theorem \ref{cb}.
\erem

\subsection{Integrable operator-valued functions} Keeping some of the notation of the previous example, let $L_1 = L_1(\RR^d, \cB)$ be the algebra of Bochner integrable $\cB$-valued functions with respect to convolution. Adjoining a formal unit $e$ to this algebra we get a unital Banach algebra denoted by $\widetilde{L_1} = \widetilde{L_1}(\RR^d, \cB)$. The algebra of Fourier transforms of elements in
$\widetilde{L_1}(\RR^d, \cB)$ is then identified with a (not closed) subalgebra  $\F(\widetilde{L_1}) \subset C_b$, with
$\hat e \equiv 1_{\cB}$ and $\hat\Phi$ vanishing at infinity for $\Phi\in L_1$. Moreover,
since $(\widehat{f*\hat\Phi})(s) = 2\pi \hat f(s)\Phi(-s)$, $s\in\RR^d$, $f\in L_1(\RR^d)$, we have
\[\|\hat\Phi\|_{1,1} = 5^d \sum_{n\in\ZZ^d} \|\phi^d_{1,n} \hat\Phi\|
\le 5^d \int_{\RR^d}\int_{[-1,1]^d}\|\Phi(s-a)\|dsda <\infty\]
and, hence, $\F(\widetilde{L_1})\subset \W(C_b)$.

The following result is a version of the celebrated Bochner-Phillips theorem \cite{BP42} for the
algebra $\widetilde{L_1}$.

\bt
Consider $A = \a e+\Phi \in \widetilde{L_1}$. Then $A$ is invertible in $\widetilde{L_1}$
if and only if $\hat A(\l) = \a1_\cB + \hat\Phi(\l)$ are invertible in $\cB$ for all $\l\in\RR^d$.
\et

\bpf
Since $\hat\Phi$ vanishes at infinity, Theorem \ref{cb} applies for $\hat A\in\W(C_b)$ and, therefore, $\hat A^{-1} \in \W(C_b)$. We need to prove that $\hat A^{-1}\in 
\F(\widetilde L_1)$.

Consider a \bai $(\varphi_n)$, $n\in\NN$, in $L_1(\RR^d)$ such that $\supp\hat\varphi_n$ is compact and $\varphi_n$ are differentiable infinitely many times, $n\in\NN$. For example,
let $\varphi\in L_1(\RR^d)$ be some infinitely many times differentiable function with $\supp\hat\varphi$ compact and
$\hat\varphi(0)=1$, and let $\varphi_n(s) = n^d\varphi(ns)$, $s\in\RR^d$, $n\in\NN$.
Then $\Phi_n = \varphi_n*\Phi\in L^1(\RR^d,\cB)$ is a sequence of entire functions that converge to $\Phi$ in $L_1$ and, moreover, the sequence $\hat A_n (\cdot) = \a 1_\cB + \hat\Phi_n(\cdot)$,
$n\in\NN$, converges to $\hat A$ in $\W(C_b)$.
Hence, there exists $N\in\NN$ such that for all $n\ge N$ elements
$\hat A_n$ are invertible in $\W(C_b)$ and $\sup_{n\ge N} \|\hat A_n^{-1}\|_{1,1}< \infty$. Since
$\supp\hat\Phi_n$ is compact and
\[(\hat A_n(\l))^{-1} = \frac1\a 1_\cB- \frac1\a\hat\Phi_n(\l)(\hat A_n(\l))^{-1}, \l\in\RR^d, n\ge N,\]
we may conclude that functions $\hat F_n :=  \frac1\a\hat\Phi_n\hat A_n^{-1}$, $n\ge N$, have compact support and therefore are Fourier transforms of some functions $F_n\in L_1$, $n\ge N$.
Moreover, by construction, the sequence $(F_n)$, $n\ge N$, is Cauchy and, therefore, converges
to some $F\in L_1$. It remains to observe that $(\hat A(\l))^{-1} = \frac1\a 1_\cB-\hat F(\l)$, $\l\in\RR^d$,
and, hence, $A^{-1} = \frac1\a e - F \in \widetilde{L_1}$.

The converse statement is trivial and the theorem is proved.
\epf

\subsection{Convolution Operators and Multipliers} 

The following example is related to the previous two but provides a view from a different vantage point.

Consider an operator valued function $\s\in L^\infty(\RR^d, B(\HH))$, where $\HH$ is a (complex) Hilbert space. Let $A_\s$ be an operator defined by
$A_\s x= \mathfrak F(\s\mathfrak F^{-1}x)$, $x\in L^2(\RR^d,\HH) = \X$, where $\mathfrak F$ is the Fourier transform operator, so that $\|\mathfrak F\|\|\mathfrak F^{-1}\|=1$. Clearly, $A_\s\in B(L^2(\RR^d,\HH))$ and $\|A_\s\| = \|\s\|_\infty$ -- the norm of the function $\s$ in $L^\infty$. Observe also that if $x_0 \in C_b(\RR^d,\HH)$ is infinitely differentiable and has compact support then $A_\s x_0 = \hat\s*x_0$ is the convolution with the Fourier transform of $\s$ in the sense of tempered
distributions.

Consider the modulation representation $\T = M$: $\RR^d\to B(\X)$ defined by \eqref{mod1} and the corresponding representation $\T_{\X\X} = \widetilde{M}$: $\RR^d\to B(B(\X))$ defined as in Example \ref{typop}, that is 
\[\widetilde{M}(t) A= M(t)AM(-t), \ t\in\RR^d, A\in B(\X). \]
Applying the above formula to the operator $A_\s$ we get $\widetilde{M}(t)A_\s x =
\mathfrak F((T(t)\s \mathfrak F^{-1} x)$, $x\in\X$, where $T$ is the translation representation defined by \eqref{trans1}. Hence, given $f\in L_1(\RR^d)$ we get $\widetilde{M}(f)A_\s x =
A_{f*\s} x$ and, therefore, 
\[\|\widetilde{M}(f)A_\s\| = \|A_{f*\s}\| = \|f*\s\|_\infty.\]
The above equalities imply that the Wiener class $\W(\X)$ contains the operator $A_\s$
if and only if $\s\in\W(C_b)$. It also follows that if $\s\in\W(C_b)$ is invertible in $L^\infty(\RR^d, B(\HH))$ then $A_\s$ is invertible and $A^{-1}_\s = A_{\s^{-1}}\in \W(\X)$.

If $\X = L^1(\RR^d, \HH)$, then the operator $A_\s$ may be defined directly as the convolution
$A_\s x = \hat\s *x$, $x\in\X$,  provided that $\hat\s\in L^1(\RR^d, B(\HH))$. In this case, $\|A_\s\|_{\W} = \|\hat\s\|_1$. This shows that if one considers convolution operators on $L^p$, Wiener class
depends non-trivially on $p\in[1,\infty]$. Recall for comparison that in case of matrices the most typical Wiener class consists of matrices with summable diagonals, i.e., it is independent of the choice of $p$ when the matrix acts on $\ell^p$. The more general approach we employed in this paper provides more flexibility and allows us to capture a greater variety of classes of memory decay.

\subsection{Integral Operators.}\label{intops} In this example we consider integral operators 
in $B(C_b(\RR^d, \CC^m))=B(\X)$. We shall impose module structure on $B(\X)$ as in Example \ref{typop}; we will write $\T_{\X\X} = \widetilde M$ if $\T_\X = M$ and
 $\T_{\X\X} = \widetilde T$ if $\T_\X = T$. 
 We begin by describing the properties of the class of operators in question. The following two propositions arise as a special case of \cite[Theorem VI.7.1]{DS88I}.

\bp\label{repint}
Assume that $A\in B(C_b)$ satisfies
\begeq\label{baicom}
\lim_\a\|M(f_\a)A - AM(f_\a)\| = 0,
\eq
where $(f_a)$ is a \bai in $L_1(\RR^d)$.
Then there exists a family $\mu_t\in Matr_{m}(M_1(\RR^d))$ of $m\times m$ matrices of finite Borel measures on $\RR^d$ such that 
\begeq\label{measrep}
(Ax)(t) = \int_{\RR^d} \mu_t(ds)x(s), x\in C_b, t\in\RR^d. \eq
\ep

\bpf
Condition \eqref{baicom}  ensures not only that $AC_0\subseteq C_0$, where $C_0$ is the space of ($\CC^m$-valued) functions vanishing at infinity, but also that $A$ is completely determined by its
action on $C_0$. Hence, \eqref{measrep} follows from
\cite[Theorem VI.7.1]{DS88I}.
\epf

\bp\label{repint1}
Assume that $A\in B(C_b)$ satisfies
\eqref{baicom}. Let $\mu_t$, $t\in\RR^d$, be the measure matrices in the representation \eqref{measrep} of $A$. Then the following are equivalent:
\begin{description}
\item[(i)]  the function $t\mapsto \mu_t: \RR^d\to Matr_{m}(M_1(\RR^d))$ is continuous;
\item[(ii)] $M(\psi)A$ is a compact operator for any $\psi\in L_1(\RR^d)$; 
\item[(iii)] $AM(\psi)$ is a compact operator for any $\psi\in L_1(\RR^d)$. 
\end{description} 
\ep




\bpf
(i) $\Leftrightarrow$ (ii) follows immediately from \cite[Theorem VI.7.1]{DS88I} in view of Proposition
\ref{repint}.

(ii) $\Rightarrow$ (iii). If $M(f_\a)A$,
where $(f_\a)$ is a b.a.i., 
are  
compact operators then, for any $\psi\in L_1(\RR^d)$,
\[AM(\psi) = \lim_\a AM(\psi*f_\a)= \lim_\a AM(\psi)M(f_\a) = \lim_\a M(f_\a)AM(\psi)\]
is a uniform limit of compact operators and, therefore, is compact. 

(iii) $\Rightarrow$ (ii) follows by  a similar argument.
\epf

\bp
Assume that $A\in B(C_b)$  satisfies
\eqref{baicom}.  Then $A$ has a representation
\begeq\label{intop1}
(Ax)(t) = \int_{\RR^d} G(t,s)x(s)ds, x\in C_b,\eq
with  $G(t,\cdot)\in Matr_{m}(L_1(\RR^d))$ for all $t\in \RR^d$, if and only if 
\[\lim_\a AM(h_a) = 0\]
 (in the uniform operator topology) for any $\g$-net $(h_\a)$, $\g\in\RR^d$.
\ep

\bpf
Let $\mu_t$, $t\in\RR^d$, be the measure matrices in the representation \eqref{measrep} of $A$. 
All these measures are absolutely continuous if and only if  
for every $\eps > 0$ there exists $\d > 0$ such that
$\sup \|Ax\| <\eps$ where the supremum is taken over all $x\in C_b$ such that $diam(\supp x) < \d$. 
The last condition is equivalent to the one we stated by Definition \ref{gnet} of $\g$-nets and the properties of the  Beurling spectrum in Lemma \ref{sprop}.
\epf

\bd
We define a class of operators $\mathfrak I\subset B(C_b)$ by requiring that each $A\in \mathfrak I$ satisfy the following conditions:
\begin{description}
\item[(i)] $A\in  (B(\X),\widetilde M)_c$;
\item[(ii)] $AM(h_\a)$ converges to $0$ in the uniform operator topology for any $\g$-net $(h_\a)$, $\g\in\RR^d$;
\item[(iii)] $M(\psi)A$ is a compact operator for each $\psi\in L_1(\RR^d))$.
\end{description} 
\ed

Observe that in view of Proposition \ref{repint1} we can replace (iii) with the equivalent requirement:
\begin{description}
\item[(iii')] $AM(\psi)$ is a compact operator for each $\psi\in L_1(\RR^d))$.
\end{description} 

We note that condition \eqref{baicom} is weaker than $A\in (B(C_b),\widetilde M)_c$ as we have shown in \cite[Section 5]{BK05}.
Hence, combining the propositions we stated in this subsection we get the following characterization of the class $\mathfrak I$.

\bp
We have $A\in\mathfrak I$ if and only if $A$ has an integral representation \eqref{intop1}, where the kernel $G$ is such that the function $t\mapsto G(t,\cdot)\!: \RR^d\to Matr_{m}(L_1(\RR^d))$ is continuous and bounded. Moreover, for such an $A$ we have
\[\|A\| = \sup_{t\in\RR^d} \int_{\RR^d} \|G(t,s)\|ds <\infty,\]
if we choose to use $|\cdot|_\infty$ as a norm on $\CC^m$.
\ep

Clearly, $\mathfrak I$ is a closed subalgebra of $B(\X)$. 
Moreover, the following lemma follows immediately from the conditions (i), (ii), (iii), and (iii').

\bl\label{rep1}
The set $\mathfrak I$ is a left ideal in $(B(\X),\widetilde M)_c$.
\el

The following key proposition is now immediate.

\bp\label{rep2}
Assume $A\in\mathfrak I$ and $I+B = (I+A)^{-1} \in B(C_b)$. Then $B\in\mathfrak I$.
\ep

\bpf
Since $I = (I+A)(I+B) = (I+B)(I+A)$, we get $B = - A - AB = -A - BA$ and $AB = BA$. Since $B$ is obviously in $(B(\X),\widetilde M)_c$, Lemma \ref{rep1} applies, and we get $BA\in\mathfrak I$. Finally, since $\mathfrak I$ is a vector space, we get $B\in\mathfrak I$.
\epf

Let us now obtain the definitions of various classes of memory decay of operators in $\mathfrak I$ with respect to the representation $\widetilde M$ in terms of the kernel $G$.   After a straightforward computation, we have
\[\W = \W(\mathfrak I, \widetilde M) = \left\{A\in\mathfrak I\!: \int_{\RR^d}\sup_{t\in\RR^d}\left(
\int_{|t-s-a|\le 1}\left\|G(t,s)\right\|ds\right)da < \infty\right\};\]
\[\bs
\W_{exp} =\! & \bigcup_{\a\in\RR^d_+}\! D(\ee) =  \Big\{A\in\mathfrak I\!: \sup_{t\in\RR^d}
\int_{|t-s-a|\le 1}\left\|G(t,s)\right\|ds
\le M e^{-\varepsilon|a|}, \\ &
\mbox{ for all } a\in \RR \mbox{ and some positive } M=M(x) \mbox{ and } \varepsilon = \varepsilon(x)
\Big\};
\end{split}
\]
\[\bs
\W_{exp}^\pm =\! & \bigcup_{\a\in\RR^d_\pm}\! D(\eet) =  \Big\{A\in\mathfrak I\!: \sup_{t\in\RR^d}
\int_{|t-s-a|\le 1}\left\|G(t,s)\right\|ds
\le M e^{-\varepsilon|a|}, \\ &
\mbox{ for all } a\in \RR_\pm \mbox{ and some positive } M=M(x) \mbox{ and } \varepsilon = \varepsilon(x)
\Big\};
\end{split}
\]
\[\B = \B(\mathfrak I, \widetilde M) = 
\left\{A\in\mathfrak I\!: \sum_{k\in\ZZ^d} \max_{|n|_\infty\ge |k|_\infty}\sup_{t\in\RR^d}
\int_{|t-s-n|\le 1}\left\|G(t,s)\right\|ds < \infty\right\}.\]
The Sobolev-type classes are, for $m\in\NN$,
\[\bs
\W^{(m)} =\ & \W(\mathfrak I, \widetilde M) = 
\Big\{A\in\mathfrak I\!:  \\
& 
\int_{\RR^d}\sup_{t\in\RR^d}\left(
\int_{|t-s-a|\le 1}\left\|G(t,s)\right\|(1+|t-s|)^mds\right)da < \infty\Big\};
\end{split}
\]
\[\bs
\B^{(m)} =\ & \B^{(m)}(\mathfrak I, \widetilde M) = 
\Big\{A\in\mathfrak I\!:  \\
& 
\sum_{k\in\ZZ^d} \max_{|n|_\infty\ge |k|_\infty}
\sup_{t\in\RR^d}
\int_{|t-s-n|\le 1}\left\|G(t,s)\right\|(1+|t-s|)^mds < \infty\Big\}.
\end{split}
\]




Next, we obtain the definitions of various classes of memory decay of operators in $\mathfrak I$ with respect to the representation $\widetilde T$. As usually, we employ
functions $\phi_{1,a}^d$ defined by \eqref{triangled}. We have
\[\bs
\W =\ & \W(\mathfrak I, \widetilde T) = \Big\{A\in\mathfrak I\!:\\
&  \int_{\RR^d}\sup_{t\in\RR^d}
\int_{\RR^d}\left\|\int_{\RR^d} \phi_{1,a}^d(u)G(t-u,s-u)du\right\|dsda < \infty\Big\};
\end{split}\]
\[\bs
\W_{exp} =\ &\W_{exp}(\mathfrak I, \widetilde T)= \bigcup_{\a\in\RR^d_+}\! D(\ee)  \\
=\ & \Big\{A\in\mathfrak I\!: \sup_{t\in\RR^d}
\int_{\RR^d}\left\|\int_{\RR^d} \phi_{1,a}^d(u)G(t-u,s-u)du\right\| ds
\le M e^{-\varepsilon|a|}, \\  &
\mbox{ for all } a\in \RR \mbox{ and some positive }  M=M(x) \mbox{ and } \varepsilon = \varepsilon(x)
\Big\};
\end{split}
\]
\[\bs
\W_{exp}^\pm =\ &\W_{exp}^\pm(\mathfrak I, \widetilde T)= \bigcup_{\a\in\RR^d_\pm}\! D(\eet)  \\
=\ & \Big\{A\in\mathfrak I\!: \sup_{t\in\RR^d}
\int_{\RR^d}\left\|\int_{\RR^d} \phi_{1,a}^d(u)G(t-u,s-u)du\right\| ds
\le M e^{-\varepsilon|a|}, \\  &
\mbox{ for all } a\in \RR_\pm \mbox{ and some positive }  M=M(x) \mbox{ and } \varepsilon = \varepsilon(x)
\Big\};
\end{split}
\]
\[\bs
 \B =\ &  \B(\mathfrak I, \widetilde T) = \Big\{A\in\mathfrak I\!: \\
&
 \sum_{k\in\ZZ^d} \max_{|n|_\infty\ge |k|_\infty}\sup_{t\in\RR^d}
\int_{\RR^d}\left\|\int_{\RR^d} \phi_{1,n}^d(u)G(t-u,s-u)du\right\|ds < \infty\Big\};
\end{split}
\]
\[\bs
\W^{(m)} =\ & \W^{(m)}(\mathfrak I, \widetilde T) = 
\Big\{A\in\mathfrak I\!:  \\
& 
\int_{\RR^d}\sup_{t\in\RR^d}
\int_{\RR^d}\left\|\int_{\RR^d} \phi_{1,a}^d(u)\frac{d^m}{du^m}G(t-u,s-u)du\right\|dsda < \infty\Big\};
\end{split}
\]
\[\bs
& \B^{(m)} =  \B^{(m)}(\mathfrak I, \widetilde T) = 
\Big\{A\in\mathfrak I\!:  \\
&
\sum_{k\in\ZZ^d} \max_{|n|_\infty\ge |k|_\infty}
\sup_{t\in\RR^d}
\int_{\RR^d}\left\|\int_{\RR^d} \phi_{1,n}^d(u)\frac{d^m}{du^m}G(t-u,s-u)du\right\|dsda < \infty\Big\}.
\end{split}
\]

In what follows we denote by $\F(\mathfrak I) = \F(\mathfrak I, \T_{\X\X})$ one of the above classes (with respect to
either $\T_{\X\X} =\widetilde M$ or $\T_{\X\X} =\widetilde T$).

\bt\label{intop}
Assume that $A \in \F(\mathfrak I)$ is such that $I+A$ is invertible in $B(C_b)$. Then $B = (I+A)^{-1} - I\in \F(\mathfrak I)$.
\et

\bpf
The result follows immediately from Theorems \ref{1sidexpwin}, \ref{invexp}, \ref{maint}, \ref{absB}, \ref{sobolt}, and Proposition \ref{rep2}.
\epf

\brem
We can easily define the classes $\F(\mathfrak I, W)$ where $W = \T_{\X\X}$ is the Weyl representation
defined as in \eqref{Weyl}. This way, however, we will not get anything significantly different from
Theorem \ref{intop}  because $\F(\mathfrak I, W) = \F(\mathfrak I, \widetilde M) \cap
\F(\mathfrak I, \widetilde T)$.
\erem

\bex
An operator $A\in \mathfrak I$ is called $2\pi$-\emph{periodic} if $\widetilde{T}(2\pi)A=A$.
For such operators we have
\[G(t+2\pi, s+2\pi) = G(t,s),\ t,s\in\RR,\]
and the function $A_{\widetilde T}$ given by $t\mapsto \widetilde{T}(t)A: \RR\to \mathfrak I$ is $2\pi$-periodic in continuous in the uniform operator topology. Let us consider the Fourier series
\cite{BK10, B97Sib, B97Izv} of this function:
\[A_{\widetilde T}(\tau)\simeq \sum_{n\in\ZZ} A_ne^{in\tau},\ \tau\in\RR,\]
where the Fourier coefficients are given by the standard formula
\[A_n = \frac1{2\pi}\int_0^{2\pi} A_{\widetilde T}(\tau) e^{in\tau}d\tau \in\mathfrak I, \ n\in\ZZ.\]
Observe that $\widetilde T(t)A_n = e^{int}A_n$, that is the Fourier coefficients are eigenvectors of the representation $\widetilde T$. Hence, $\widetilde T(t)[A_nM(-n)] = A_nM(-n)$, $t\in\RR$,  $n\in\ZZ$, and, therefore,
\[(A_nx)(t) = \int_\RR G_n(t-s)e^{ins}x(s)ds, t\in\RR, n\in\ZZ.\]
Observe now that for $f\in L_1(\RR)$ we have $\widetilde M(f)A \simeq \sum_{n\in\ZZ} \hat f(n) A_n$.
Therefore, a $2\pi$-periodic operator $A\in \mathfrak I$ belongs to $\W(\mathfrak I,\widetilde M)$ if and only if it has summable Fourier coefficients, i.e. 
\[\sum_{n\in\ZZ} \|A_n\| = \sum_{n\in\ZZ} \|G_n\|_{L_1(\RR)} < \infty.\]
\eex

\brem
We refer to \cite{FS10, K99, K01, SS09, S08} for related results involving  different classes of integral operators.
\erem

\subsection{Integro-differential operators with delay} Here we consider operators $A = \frac d{dt} + 
L+C$, where 
\[(Lx)(t) =  \sum_{k=1}^\infty (L_kT(h_k) x)(t) = \sum_{k=1}^\infty \ell_k(t)x(t+h_k),\]
$ x\in C_b(\RR,\CC^m)$,   $k\in\NN$, $h_k \in \RR$,  $\ell_k\in C_b(\RR,B(\CC^m))$, $\sum\limits_{k\in\NN} L_kT(h_k)$ converges in the uniform operator topology, and
\[(Cx)(t) = \int_\RR G(t,s)x(s)ds,\] 
is an integral operator from Subsection \ref{intops}, i.e., $C\in \mathfrak I$. We will consider
$A$ as a bounded operator from the space $\X = C^1_b(\RR, \CC^m)$ of bounded continuously differentiable functions into $\Y = C_b(\RR, \CC^m)$ as well as an unbounded operator 
$A: \X =D(A)\subset \Y\to\Y$.

\bl
Assume that $B = A^{-1} \in B(\Y)$. Then $B\in\mathfrak I$.
\el

\bpf
Let $D=\frac d{dt} + I$: $D(A)\subset\Y\to\Y$. Clearly, $D^{-1} \in B(\Y)$ and
\begeq\label{din}
(D^{-1}x)(t) = \int_\infty^t e^{s-t}x(s)ds,
\eq
so that $D^{-1}\in\mathfrak I$. Hence, we have
\begeq\label{nibd}
B = (L+C+D - I)^{-1} = D^{-1}(I +(L+C-I)D^{-1})^{-1} \in \mathfrak I\eq
by Lemma \ref{rep1} and Proposition \ref{rep2} which are both applicable because convergence of $\sum\limits_{k\in\NN} L_kT(h_k)$ in the uniform operator topology implies $L\in (B(\Y), \widetilde M)_c$.
\epf

Firstly, we consider the module structure on $L(\X,\Y)$ associated with the representation $\T_{\X\Y}$ defined as in the example \ref{typop} where $\T_\X$ and $\T_\Y$ are the
translation representations given by \eqref{trans1}. With a slight abuse of notation we will  use the symbol $\widetilde T$ to denote either $\T_{\X\Y}$, $\T_{\Y\X}$, $\T_{\X\X}$ or $\T_{\Y\Y}$. It will be clear from the context which of the representations may be used in each instance. 
Let $\F$ denote one of the classes of operators with memory decay considered in this paper. 
The following
 conditions are sufficient to ensure $A\in \F = \F(L(\X,\Y),\widetilde T)$: 
\begin{description}
\item[(A)] $L_k \in \F(C_b, T)$ for each $k\in\NN$, and $\sum\limits_{k\in\NN}\|L_k\|_{\F} < \infty$;
\item[(B)] $C\in\F(\mathfrak I, \widetilde T)$.
\end{description}

\bt
Assume that an operator $A = \frac d{dt} + L+C$ satisfies conditions \emph{(A)} and \emph{(B)} and
$B = A^{-1} \in B(\Y)$. Then $B \in \F(\mathfrak I,\widetilde T)$. Moreover, if we regard $B\in L(\Y,\X)$ then
$B\in \F(L(\Y,\X),\widetilde T)$.
\et

\bpf
Observe that \eqref{din} implies that $D^{-1}\in \F$ for any $\F$. Hence, the
result follows from \eqref{nibd} and Theorem 
\ref{intop}.
\epf

Secondly, we consider the module structure on $L(\X,\Y)$ associated with the representation $\T_{\X\Y}$ defined as in the example \ref{typop} where $\T_\X$ and $\T_\Y$ are the
translation representations given by \eqref{mod1}. We will  use the symbol $\widetilde M$ to denote either $\T_{\X\Y}$ 
or $\T_{\Y\Y}$ with the understanding that the context will determine which of the representations may be used in each instance. 
It is not hard to determine  when 
 $L+C \in \F = \F(B(\Y),\widetilde M)$.
 In all cases we would need $C\in\F(\mathfrak I, \widetilde M)$.  Conditions on the operator $L$ vary:
\[ L\in \W(B(\Y),\widetilde M)  \mbox{ if and only if } 
\sum\limits_{n\in\ZZ}\left\|\sum\limits_{k\in\ZZ} \hat\phi^d_{1,n}(h_k)\ell_k\right\| <\infty;\]
\[L\in \W_{exp}(B(\Y),\widetilde M) \mbox{ if and only if }  
\left\|\sum\limits_{k\in\ZZ} \hat\phi^d_{1,n}(h_k)\ell_k\right\| \le M\g^{|n|}, n\in\NN,\]
for some $M > 0$, and $\g\in(0,1)$,  
and so on; we leave the remaining classes to the reader.

\bt
Assume that an operator $A = \frac d{dt} + L+C$ satisfies $L\in \F(B(\Y),\widetilde M)$ and
$C\in\F(\mathfrak I, \widetilde M)$. If $B = A^{-1}\in  B(\Y)$, then $B \in \F(\mathfrak I,\widetilde M)$. 
\et

\bpf
The
result follows from \eqref{din}, \eqref{nibd} and Theorem 
\ref{intop}.
\epf

\brem
Observe that unlike $\T_\Y$, $\T_\X$ is not a bounded representation in this case (see \eqref{mod1}). Hence, we cannot simply quote the general results of this paper to obtain the memory decay for $B = A^{-1}\in L(\Y,\X)$. First, one would have to modify the
general definitions of the classes $\F$ considered in this paper. 
For an operator given by \eqref{intop1}
description of the classes $\F(L(\Y,\X),\widetilde T)$ and $\F(L(\Y,\X),\widetilde M)$ is 
obtained using the formulas for $\F(\mathfrak I)$ with the kernel $G$ replaced by 
%
$DG = G+\frac{\partial G}{\partial t}$. For more general operators more work would be needed.
Secondly, one would have to prove an analog of Lemma 
\ref{algprop} and similar results for other classes of memory decay. We will pursue these and other extensions in a sequel to this paper.
\erem
\medskip
\centerline{\bf Acknowledgements}
\medskip
We are grateful to K. Gr\"ochenig and Q.~Sun for inspiration provided by their papers and  useful comments that helped us improve the presentation.

\bibliographystyle{siam}
\bibliography{../refs}

\end {document}